\documentclass[preprint, review, 12pt]{elsarticle}

\usepackage{amssymb}
\usepackage{fullpage}
\usepackage{amsthm}
\usepackage{amsmath}
\usepackage{caption}
\usepackage{subcaption}
\usepackage{float}
\usepackage{hyperref}
\usepackage{bm}
\usepackage{float}
\usepackage{xcolor}
\usepackage{physics}
\usepackage{soul}
\usepackage{mathtools}

\journal{Journal Name}

\begin{document}

\begin{frontmatter}

\title{Simulation of non-stationary and non-Gaussian random processes by 3rd-order Spectral Representation Method: Theory and POD implementation}

\author{Lohit Vandanapu}
\author{Michael D. Shields}

\address{Johns Hopkins University, Baltimore, United States}

\begin{abstract}

This paper introduces the $3^{rd}$-order Spectral Representation Method for simulation of non-stationary and non-Gaussian stochastic processes. The proposed method extends the classical $2^{nd}$-order Spectral Representation Method to expand the stochastic process from an evolutionary bispectrum and an evolutionary power spectrum, thus matching the process completely up to third-order. A Proper Orthogonal Decomposition (POD) approach is further proposed to enable an efficient FFT-based implementation that reduces computational cost significantly. Two examples are presented, including the simulation of a fully non-stationary seismic ground motion process, highlighting the accuracy and efficacy of the proposed method.

\end{abstract}

\begin{keyword}
Spectral Representation \sep Non Stationary \sep Fast Fourier Transform
\end{keyword}

\end{frontmatter}

\section{Introduction}
\label{sec:Introduction}

Monte Carlo Simulation of non-stationary stochastic processes is of extreme importance, especially for simulation of extreme events like earthquake ground motion or transient wind gusts. These simulations are particularly essential in the case of structural systems involving non-linear dynamics where uncertainity cannot be quantified analytically.

Properties of a stationary process remain the same at every instant of time and because to this independence in frequency and time, computationally efficient methods for the simulation of stationary process can be derived. On the other hand the properties of a non-stationary process vary with time making their simulation challenging because this introduces dependence between time and frequency. For the simulation of stationary stochastic processes, the Spectral Representation Method (SRM) and Karhunen-Loeve Expansion (KLE) \cite{Huang2001} are the most widely-used methods. Both methods rely on an expansion of a random process having a truncated form:
\begin{equation}
    X(t) \approx \sum_{i=1}^{N} \theta_{i}(\omega) \psi_{i}(t)
\label{eqn:general_expansion}
\end{equation}
where $\psi(t)$ are orthogonal basis functions and $\theta(\omega), \omega \in \Omega$ are random variables. The primary difference between SRM and KLE is that SRM uses harmonic basis functions, whereas KLE uses the eigen functions of the covariance function as the basis functions. In this work, we will focus on novel developments for the SRM.

Shinozuka was the first to use the Spectral Representation Method for the simulation of stochastic processes \cite{Shinozuka1972}. Then, Yang \cite{Yang1972} suggested the use of Fast Fourier Transform (FFT) to improve its computational efficiency. The complete theoretical footing for the SRM was then established in the 1990s along with extensions for the simulation of ergodic vector processes\cite{Deodatis1996a}, stochastic waves \cite{Deodatis1989}, random fields \cite{Shinozuka1996}.

Samples generated by the SRM are Gaussian (either naturally or asymptotically depending on the implementation \cite{grigoriu1993spectral}). This is convenient for many applications, but limiting for applications involving non-Gaussian (i.e. higher-order) stochastic processes. To fill this gap, several methods for the simulation of non-Gaussian processes have been proposed. One class of such methods works by introducing correlated random variable with deterministic basis functions like Hermite and Legendre polynomials \cite{Puig2002, Liu2017}. Another class of processes works by performing an inverse Cumulative Distribution Function (CDF) transform of Gaussian samples \cite{Grigoriu1998} referred to as the 'translation' process as
\begin{equation}
    Y(t) = F^{-1}(\Phi(X(t)))
\end{equation}
where $X(t)$ is a standard Gaussian stochastic process, $\Phi(\cdot)$ is the standard normal CDF and $F(\cdot)$ is the CDF of the non-Gaussian distribution. A wide range of methods utilising translation process theory have been developed \cite{Shields2011, Shields2013, Kim2015}. More recently, a class of methods has been developed that theoretically extends the SRM to higher-order ($3^{rd}$-order asymmetrically non-Gaussian) stochastic processes by considering interactions among the wave components of the SRM expansion \cite{Shields2017, Vandanapu2021}.

The SRM for simulation of non-stationary stochastic processes was also presented in the seminal paper by Shinozuka \cite{Shinozuka1972} and relies on the theory of evolutionary power developed by Priestley \cite{Priestley1965}. But the formula was based on the summation of the trigonometric cosines functions which is computationally inefficient. To improve the efficiency of the simulation process, Li and Kareem \cite{Li1997} proposed the use of discrete fourier transform in tandem with digital filtering. Huang \cite{Huang2014} also developed FFT-aided methods involving the use of wavelets. All these methods focused only on using FFT directly on the evolutionary power spectrum. But, additional savings can be obtained by projecting the evolutionary power spectrum onto a lower-dimensional orthogonal basis in the frequency domain and considering time varying modulations of the corresponding components. This has the effect of decomposing the non-stationary process into the sum of time-frequency separable processes in which each process has a stationary component that is modulated in time.  Li and Kareem \cite{Li1991} suggested the use of orthogonal legendre polynomials to decompose the evolutionary spectrum and established an FFT-aided spectral representation method. Huang \cite{Huang2015} proposed an FFT-based approach that leverages the Proper Orthogonal Decomposition (POD) for the simulation of multivariate processes. POD was also used in conjunction with Stochastic wave theory in the simulation of multivariate non-stationary random processes \cite{Peng2017}. The POD involves projecting high-dimensional data (i.e. the evolutionary spectrum) onto to a low-dimensional manifold by finding the best set of basis functions. Since most of the information is contained within a few modes, the POD can drastically reduce the dimensionality of the data. For SRM-based simulations this drastically reduces computation cost.

Simulation of non-stationary and non-Gaussian random processes compounds the challenges of simulating non-Gaussian processes with those of simulating non-stationary processes. As a result, few methods have been developed to successfully simulate these complex processes \cite{sakamoto2002simulation,ferrante2005stochastic,shields2013estimation, Kim2015, dai2019explicit, montoya2019simulation, zheng2021sample}. Those that use the SRM for simulation, rely on an expansion from the evolutionary spectrum coupled with translation process theory \cite{ferrante2005stochastic, shields2013estimation,Kim2015}. Here, we derive a direct third-order SRM from the evolutionary power spectrum and the evolutionary bispectrum for the simulation of non-Gaussian stochastic processes having known second and third-order properties. This method extends the third-order SRM previously developed for stationary random processes \cite{Shields2017}, multi-dimensional random fields, and stochastic vector processes \cite{Vandanapu2021}. We further develop a POD-based implementation for the simulations, which allows the use of FFT and drastic computational improvement of the simulation formula. The theoretical properties of the expansion are derived and the effectiveness of the proposed methodology is demonstrated with the use of two numerical examples: one illustrating a simple time-frequency separable third-order process and one considering a non-Gaussian stochastic ground motion process.
 


\section{2nd-order spectral representation for non-stationary stochastic processes}


In general, a one-dimensional, uni-variate, zero-mean non-stationary process, $X(t)$, can be expressed as
\begin{equation}
    X(t) = \int_{-\infty}^{\infty} \phi(t, \omega) dZ(\omega)
\label{eqn:priestley_representation}
\end{equation}
where $Z(\omega)$ is a spectral process with orthogonal increments having the following properties
\begin{equation}
\begin{aligned}
    & \mathbb{E}[dz(\omega)] = 0 \\
    & \mathbb{E}[|dz(\omega)|^2] = d\mu( \omega) \\
\end{aligned}
\end{equation}
where $\phi(t, \omega)$ is selected from a suitable family of functions and $\mu(\omega)$ an associated measure such that the covariance function of the process can be expressed as:
\begin{equation}
    Cov(X(t),X(s))=E[X(t)X(s)]=\int_{-\infty}^\infty \phi(t,\omega)\phi^*(s,\omega)d\mu(\omega).
    \label{eqn:covariance}
\end{equation}
The expression of the stochastic process in Eq. \eqref{eqn:priestley_representation} is not unique owing to the fact that $\phi(t,\omega), \mu(\omega)$ can be selected arbitrarily to satisfy Eq.\ \eqref{eqn:covariance}. A common selection for stationary stochastic processes is the complex exponentials $\phi(t,\omega)=e^{\iota\omega t}$ such that Eq.\ \eqref{eqn:priestley_representation} admits the classical Cramer spectral representation \cite{cramer1967}. The complex exponentials cannot be used for non-stationary processes. Instead, Priestley \cite{Priestley1965,Priestley1967} suggested to use the amplitude modulated complex exponentials such that
\begin{equation}
    \phi(t,\omega)=A(t,\omega)e^{\iota\omega t},
\end{equation}
referred to as an oscillatory process. With this representation, the process is expressed as
\begin{equation}
    X(t) = \int_{-\infty}^{\infty} A(t, \omega) e^{\iota \omega t} dz(\omega)
\label{eqn:priestley_representation}
\end{equation}
and the evolutionary power spectrum can be defined as
\begin{equation}
    dS(t,\omega)=|A(t,\omega)|^2d\mu(\omega)
\end{equation}
Priestley further suggests that it is convenient to standardize the modulating function such that $A(0,\omega)=1$, implying that the measure $\mu(\omega)=S(\omega)$ is equal to the power spectrum at time $t=0$ and $A(t,\omega)$ represents the time change from this original power spectrum. Under these conditions, the two-sided evolutionary power spectral density function is defined as
\begin{equation}
    S(t, \omega) = |A(t, \omega)|^2 S(\omega)
    \label{eqn:EPSDF}
\end{equation}

Given the spectral representation in Eq.\ \eqref{eqn:priestley_representation}, Shinozuka \cite{Shinozuka1972} showed that the non-stationary process can be simulated by 
\begin{equation}
    X(t) = \sqrt{2} \sum_{n=0}^{N-1} \sqrt{2S(t, \omega_n)\Delta \omega} \cos(\omega_n t + \Phi_{n})
\label{eqn:2_order_esrm}
\end{equation}
where $\Delta \omega$ is the frequency interval with
\begin{equation}
\begin{aligned}
    & \omega_n = n\Delta\omega \\
    & \Delta\omega = \frac{\omega_u}{N}. \\
\end{aligned}
\end{equation}
and $\omega_{u}$ is the cutoff frequency. With the evolutionary power spectral density function given in Eq.\ \eqref{eqn:EPSDF}, this simulation equation can be equivalently expressed as:
\begin{equation}
    X(t) = \sqrt{2} \sum_{n=0}^{N-1} A(t,\omega_n)\sqrt{2S( \omega_n)\Delta \omega} \cos(\omega_n t + \Phi_{n})
\label{eqn:2_order_esrm}
\end{equation}

\section{3rd-order spectral representation for non-stationary stochastic processes}

In this section, we extend the spectral representation theory to third-order non-stationary stochastic processes. We then present a new expression for the simulation of third-order non-stationary stochastic processes that leverages the spectral representation.

\subsection{3rd-order non-stationary spectral representation}

The non-stationary spectral representation in Eq.\ \eqref{eqn:priestley_representation} can, in general, be extended to represent stochastic processes of arbitrary order by extending the orthogonality conditions on the spectral process $Z(\omega)$. Extension to third-order introduces the following orthogonality conditions
\begin{equation}
\begin{aligned}
    & \mathbb{E}[dZ(\omega)] = 0 \\
    & \mathbb{E}[|dZ(\omega)|^2] = S(\omega) d \omega \\
    & \mathbb{E}[dZ(\omega_1)dZ(\omega_2)dZ^{*}(\omega_3)] = \delta(\omega_{1} + \omega_{2} - \omega_{3}) B(\omega_{1}, \omega_{2}) d\omega_1 d\omega_2 
\label{eqn:3_order_orthogonal}
\end{aligned}
\end{equation}
Notice that the measure $d\mu(\omega)=S(\omega)d\omega$ in the second-order condition and that the third-order condition introduces an associated stationary bispectrum $B(\omega_1,\omega_2)$. Analogous to Eq.\ \eqref{eqn:EPSDF}, we define the evolutionary bispectrum as
\begin{equation}
\begin{aligned}
    B(t, \omega_1, \omega_2) &= A(t, \omega_1)A(t, \omega_2)A(t, \omega_1 + \omega_2) \mathbb{E}[dZ(\omega_1)dZ(\omega_2)dZ^*(\omega_1 + \omega_2)]\\
    B(t, \omega_1, \omega_2) &= A(t, \omega_1)A(t, \omega_2)A(t, \omega_1 + \omega_2) B(\omega_1,\omega_2)
\end{aligned}
\end{equation}
such that it represents the distribution of the skewness over the space of frequency pairs at any given time, analogous to Priestley's \cite{Priestley1965} evolutionary spectrum distributing the variance over the frequency domain at any given time. This is shown in more detail in \ref{ap:C}.







Next, consider that the orthogonal increments $dZ(\omega)$ can be divided into their real and complex components as
\begin{equation}
\begin{aligned}
    dU(\omega) &= \Re[dZ(\omega)],\\ 
    dV(\omega) &= -\Im [dZ(\omega)]
\end{aligned}
\end{equation}
where the individual increments $dU(\omega)$ and $dV(\omega)$ satisfy
\begin{equation}
    \mathbb{E}[dU(\omega)] = \mathbb{E}[dV(\omega)] = 0
\end{equation}
Similarly, the modulating function $A(t, \omega)$ can be divided into real and complex components as
\begin{equation}
    A(t, \omega) = \alpha(t, \omega) + \iota \beta(t, \omega)
\end{equation}


Applying these relations, the spectral representation in \eqref{eqn:priestley_representation} can be rewritten as
\begin{equation}
    X(t) = \int_{-\infty}^{\infty}[\alpha(t, \omega) + \iota\beta(t, \omega)][\cos(\omega t) + \iota \sin(\omega t)][dU(\omega) - \iota dV(\omega)]
\end{equation}
and is equivalently represented as
\begin{equation}
    X(t) = \int_{-\infty}^{\infty}\cos(\omega t)dU_{t}(\omega) - \sin(\omega t)dV_{t}(\omega)
\label{eqn:process_in_terms_of_orthogonal_increments}
\end{equation}
where
\begin{equation}
    dU_t(\omega) = \alpha(t, \omega)dU(\omega) + \beta(t, \omega)dV(\omega)
    \label{eqn:ortho_U}
\end{equation}
\begin{equation}
    dV_t(\omega) = \beta(t, \omega)dU(\omega) - \alpha(t, \omega)dV(\omega)
    \label{eqn:ortho_V}
\end{equation}
are modulated orthogonal increments satisfying the following properties (see \ref{sec:AppendixA})
\begin{equation}
\begin{aligned}
    & \mathbb{E}[dU_{t}(\omega)] = \mathbb{E}[dV_{t}(\omega)] = 0 \\
    & \mathbb{E}[dU^{2}_{t}(\omega)] = \mathbb{E}[dV^{2}_{t}(\omega)] = 2S(\omega, t) d\omega \\
    & \mathbb{E}[dU_{t}(\omega_{1})dU_{t}(\omega_{2})dU_{t}(\omega_{1} + \omega_{2})] = \mathbb{E}[dV_{t}(\omega_{1})dV_{t}(\omega_{2})dV_{t}(\omega_{1} + \omega_{2})] = \\
    & = 2B(t, \omega_{1}, \omega_{2})d\omega_{1}d\omega_{2} \\
\label{eqn:modulated_orthogonal_properties}
\end{aligned}
\end{equation}

\subsection{Simulation of 3rd-order non-stationary stochastic processes}

To enable simulation, we propose the following modulated orthogonal increments that satisfy the required properties
\begin{equation}
\begin{aligned}
    dU_{t}(\omega_{k}) &= [2S_{p}(t, \omega_{k})\Delta \omega]^{\frac{1}{2}}\cos(\phi_{k})\\
    &+ \sum_{i+j=k}^{i \geq j \geq 0}[2S(t, \omega_{k})\Delta \omega]^{\frac{1}{2}}|b_{p}(t, \omega_{i}, \omega_{j})|\cos(\phi_{i} + \phi_{j} + \beta(t, \omega_{i}, \omega_{j}))\\
\end{aligned}
\end{equation}
\begin{equation}
\begin{aligned}
    dV_{t}(\omega_{k}) &= [2S_{p}(t, \omega_{k})\Delta \omega]^{\frac{1}{2}}\sin(\phi_{k})\\
    &+ \sum_{i+j=k}^{i \geq j \geq 0}[2S(t, \omega_{k})\Delta \omega]^{\frac{1}{2}}|b_{p}(t, \omega_{i}, \omega_{j})|\sin(\phi_{i} + \phi_{j} + \beta(t, \omega_{i}, \omega_{j}))\\
\end{aligned}
\end{equation}
where $S(t, \omega)$ is the evolutionary power spectrum, $S_p(t, \omega)$ is the pure component of the evolutionary power spectrum defined as
\begin{equation}
    S_p(t, \omega_{k}) = S(t, \omega_{k})(1 - \sum_{i+j=k}^{i \geq j \geq 0}b_{p}^{2}(t, \omega_{i}, \omega_{j}))
\label{eqn:pure_eps_definition}
\end{equation}
and $b_{p}(t, \omega_{i}, \omega_{j})$ is the partial evolutionary bicoherence given by
\begin{equation}
    b_{p}^{2}(t, \omega_{i}, \omega_{j}) = \frac{|B(t, \omega_{i}, \omega_{j})|^{2}\Delta\omega}{S_{p}(t, \omega_{i})S_{p}(t, \omega_{j})S(t, \omega_{i} + \omega_{j})}
\end{equation}
where $B(t, \omega_1, \omega_2)$ is the evolutionary bispectrum. Note that the pure evolutionary power spectrum and the evolutionary bicoherence result from a direct extension of their stationary counterparts introduced in \cite{Shields2017}.

Using the above definitions for the modulated orthogonal increments, the stochastic process in Eq.\ \eqref{eqn:process_in_terms_of_orthogonal_increments} process can be expanded as
\begin{equation}
\begin{aligned}
    X(t) &= 2\sum_{k=-\infty}^{\infty}\sqrt{S_{p}(t, \omega_{k})\Delta\omega}\cos(\omega_{k}t + \phi_{k})\\
    & + 2\sum_{k=0}^{N-1}\sum_{i+j=k}^{i \geq j \geq 0}\sqrt{S(t,\omega_{k})\Delta\omega}|b_{p}(t, \omega_{i}, \omega_{j})|\cos(\omega_{k}t + \phi_{i} + \phi_{i} + \beta(t, \omega_{i}, \omega_{j}))
\label{eqn:process_infinite_sum}
\end{aligned}
\end{equation}
where $\beta(t,\omega_i,\omega_j)$ is the evolutionary biphase given by
\begin{equation}
    \beta(t,\omega_i,\omega_j) = \arctan(\frac{\Im{B(t, \omega_1, \omega_2)}}{\Re{B(t, \omega_1, \omega_2)}})
\end{equation}
where $\Im{\cdot}$ represents the imaginary part, $\Re{\cdot}$ represents the real part and
\begin{equation}
\begin{aligned}
    \Delta \omega &= \frac{\omega_{u}}{N}\\
    \omega_{k} &= k\Delta \omega \: \: k = 1,2,\dots,N-1
\end{aligned}
\end{equation}
and $\omega_{u}$ is the upper cutoff frequency beyond which the evolutionary power spectral density function $S(t, \omega)$ may be assumed to be zero. Finally, it is assumed that $S(t,\omega_0)=0$.

It is shown in \ref{ap:B} that the above non-stationary process satisfies the correct ensemble statistical properties up to third-order.

The infinte series representation in Eq.\  \eqref{eqn:process_infinite_sum} can be truncated with to include $N$ terms for simulation purposes as 
\begin{equation}
\begin{aligned}
    X(t) &= 2\sum_{k=0}^{N-1}\sqrt{S_{p}(t, \omega_{k})\Delta\omega}\cos(\omega_{k}t + \phi_{k})\\
    & + 2\sum_{k=0}^{N-1}\sum_{i+j=k}^{i \geq j \geq 0}\sqrt{S(t,\omega_{k})\Delta\omega}|b_{p}(t, \omega_{i}, \omega_{j})|\cos(\omega_{k}t + \phi_{i} + \phi_{i} + \beta(t, \omega_{i}, \omega_{j}))
\end{aligned}
\end{equation}





Finally, after some rearrangement, the simulation formula can be expressed as
\begin{equation}
\begin{aligned}
    X(t) &= 2\sum_{k=0}^{N-1}\sqrt{S(t, \omega_{k})\Delta\omega}\Big[\sqrt{(1 - \sum_{i+j=k}^{i \geq j \geq 0}b_{p}^{2}(t, \omega_{i}, \omega_{j}))}\cos(\omega_{k}t + \phi_{k})\\
    & + \sum_{i+j=k}^{i \geq j \geq 0}|b_{p}(t, \omega_{i}, \omega_{j})|\cos(\omega_{k}t + \phi_{i} + \phi_{i} + \beta(t, \omega_{i}, \omega_{i}))\Big]
\label{eqn:sum_of_cosines_formula}
\end{aligned}
\end{equation}

The computational complexity of the simulation formula is of the order $O(MN^2)$, where $M$ is the number of time discretizations, which leads to high computational expense for simulations in this form. Given the computational expense, we propose an alternative formulation that leverages the proper orthogonal decomposition (POD) next.  


\section{POD based implementation of 3-order Spectral Representation Method}

Li and Kareem \cite{Li1991} first proposed the use of the POD technique to enable the use of the Fast Fourier Transform (FFT) for the simulation of non-stationary random processes. POD for the $2^{nd}$ Spectral Representation Method involves the decomposition of $\sqrt{S(t, \omega)}$ in Eq.\ \eqref{eqn:2_order_esrm} into a sum of separable frequency and time functions as
\begin{equation}
    \sqrt{S(t, \omega)} = \sum_{q=1}^{N_q} a_{q}(t) \Phi_{q}(\omega)
\label{eqn:pod_sqrt_ps}
\end{equation}
where $\Phi_{q}(\omega)$ are orthogonal functions and $a_q(t)$ are time-dependent principal coordinate vectors calculated by $a_q(t) = \int_{\omega} \sqrt{S(t, \omega)}\Phi_q(\omega)$. After decomposing the evolutionary power spectrum into the set of time and frequency functions $a_q(t), \Phi_q(\omega)$, the simulation formula in Eq.\ \eqref{eqn:2_order_esrm} can be expressed as
\begin{equation}
    X(t) = \sum_{q=1}^{N_q} 2 a_q(t) \sum_{n=0}^{N-1} \Phi_q(\omega_n) \sqrt{\Delta \omega} \cos(\omega_n t + \phi_{n}^{q})
\end{equation}
This is equivalent to a sum of modulated stationary random processes. Since an individual stationary random process can be simulated using FFT, the non-stationary process can be simulated as a sum of processes simulated using FFT.


In the remainder of the section, we develop the POD-based implementation for simulation of non-stationary processes by $3^{rd}$-order Spectral Representation Method. 

Recall the simulation formula in Eq.\ \eqref{eqn:process_infinite_sum}. Let us apply the Tucker decomposition \cite{Kolda2009} on the second term as
\begin{equation}
    \frac{|B(t, \omega_1, \omega_2)|}{\sqrt{S_{p}(t, \omega_1)S_{p}(t, \omega_2)}} = \mathcal{T} \times \prescript{}{1}{U^{(1)}} \times \prescript{}{2}{U^{(2)}} \times \prescript{}{3}{U^{(3)}}
    \label{eqn:tucker_decomposition}
\end{equation}
where $\mathcal{T}$ is the core tensor and $U^{(1)}, U^{(2)}, U^{(3)}$ are unitary matrices. A set of orthogonal functions $\Phi_{q}(\omega)$ is chosen to be the columns of the matrix $U^{(3)}$. From the set of basis functions $\Phi_{q}(\omega)$, we can form a second-order basis as
\begin{equation}
    \theta_{rs}(\omega_1, \omega_2) = \Phi_{r}(\omega_1) \Phi_{s}(\omega_2)
\end{equation}
This allow us to decouple the time and frequency components in the interactive component of the simulation formula, $\frac{B(t, \omega_1, \omega_2)}{\sqrt{S_{p}(t, \omega_1)S_{p}(t, \omega_2)}}$, which depend on frequency pairs. The corresponding set of projected bispectrum amplitude functions $b_{rs}(t)$ can be computed by
\begin{equation}
\begin{aligned}
    b_{rs}(t) & = \int_{\omega_1, \omega_2} \theta_{rs}(\omega_1, \omega_2) \frac{B(t, \omega_1, \omega_2)}{\sqrt{S_{p}(t, \omega_1)S_{p}(t, \omega_2)}} \\
    & = \int_{\omega_1, \omega_2} \Phi_{r}(\omega_1) \Phi_{s}(\omega_2) \frac{B(t, \omega_1, \omega_2)}{\sqrt{S_{p}(t, \omega_1)S_{p}(t, \omega_2)}} \\
\end{aligned}
\end{equation}
where the reconstruction can be expressed as
\begin{equation}
    \frac{B(t, \omega_1, \omega_2)}{\sqrt{S_{p}(t, \omega_1)S_{p}(t, \omega_2)}} \approx \sum_{r, s} b_{rs}(t) \theta_{rs}(\omega_1, \omega_2).
    \label{eqn:bispectrum_approximation}
\end{equation}
This represents the higher-order anologue to Eq.\ \eqref{eqn:pod_sqrt_ps} -- i.e. the separation of the bispectrum into orthogonal frequency functions and corresponding time-dependent modulating functions -- and further allows the pure component of the power spectrum to be decomposed as
\begin{equation}
\begin{aligned}
    & \sqrt{S_{p}(t, \omega)} \approx \sum_{q=1} a_{q}(t) \Phi_{q}(\omega) \\
    & a_q(t) = \int_{\omega} \sqrt{S_{p}(t, \omega)}\Phi_q(\omega) \\
    \label{eqn:second_order_projection}
\end{aligned}
\end{equation}

Note that the projected amplitude function $b_{rs}(t)$ function can be complex-valued since the non-stationary Bispectrum $B(t, \omega_i, \omega_j)$ can be complex-valued. Each of these functions can be expressed with a set of biphase angles $\gamma_{rs}(t)$ as
\begin{equation}
    b_{rs}(t) = |b_{rs}(t)| e^{\iota \gamma_{rs}(t)}
\end{equation}




Substituting the functions $\Phi_q(\omega), a_q(t), b_{rs}(t)$ into the simulation formula yields
\begin{equation}
\begin{aligned}
    X(t) & = 2\sum_{k=0}^{N-1} \Bigg[\sum_{q = 1}^{N_q} a_q(t)\Phi_q(\omega_k) \sqrt{\Delta \omega} \cos(\omega_k t - \phi_{kq}) \\
    & + \sum_{i+j=k}^{i \geq j \geq 0}\sum_{r=1}^{N_{q}} \sum_{s=1}^{N_{q}} |b_{rs}(t)|\theta_{rs}(\omega_i, \omega_j) \Delta \omega \cos(\omega_k t - \phi_{ri} - \phi_{sj} + \gamma_{rs}(t)) \Bigg] \\
    \label{eqn:pod_expansion}
\end{aligned}
\end{equation}
which can be further simplified to
\begin{equation}
\begin{aligned}
    X(t) & = 2 \sum_{k=0}^{N-1} \sum_{r = 1}^{N_q} \Bigg[a_r(t) \Phi_r(\omega_k) \sqrt{\Delta \omega} \cos(\omega_k t - \phi_{kr}) \\
    & + \sum_{i+j=k}^{i \geq j \geq 0} \sum_{s=1}^{N_{q}} |b_{rs}(t)|\theta_{rs}(\omega_i, \omega_j) \Delta \omega \cos(\omega_k t - \phi_{ri} - \phi_{sj} + \gamma_{rs}(t)) \Bigg] \\
    X(t) & = 2\sum_{r = 1}^{N_q} x_r(t)\\
\label{eqn:non_stationary_sum_fft}
\end{aligned}    
\end{equation}
where
\begin{equation}
\begin{aligned}
    x_r(t) & = \sum_{k=0}^{N-1} \Big[ a_r(t) \Phi_{r}(\omega_k) \sqrt{\Delta\omega} \cos(\omega_k t - \phi_{kr}) \\
    & + \sum_{i+j=k}^{i \geq j \geq 0} \sum_{s=1}^{N_q} |b_{rs}(t)| \theta_{rs}(\omega_i, \omega_j) \Delta \omega \cos(\omega_k t - \phi_{ri} - \phi_{sj} + \gamma_{rs}(t)) \Big]
\end{aligned}
\end{equation}
That is, the third-order non-stationary process can be simulated as the sum of modulated third-order stationary processes, $x_r(t)$. We observe from \cite{Vandanapu2021} that $x_{q}(t)$ can be simulated using FFT, and therefore the full non-stationary third-order stochastic process $X(t)$ can also be simulated using FFT.


The processes simulated by the POD method satisfies the required second and third properties as shown in \ref{ap:D} and \ref{ap:E}. Finally, the number of components $N_q$ in the expansion from Eq.\ \eqref{eqn:tucker_decomposition} should then be chosen such that the approximation Eq.\ \eqref{eqn:bispectrum_approximation} is sufficiently accurate. As we show in the numerical examples that follow, this number $N_q$ is often small. The computational complexity of the POD based implementation is of the order
\begin{equation}
    O(N_{q}^{2}N^{2} + N_{q}^{2}M\log{N}),
    \label{eqn:pod-time-complexity}
\end{equation}
which provides considerable time savings over the original simulation equation, which again is $O(MN^2)$, as long as $N_{q}^{2} << M$.







\section{Numerical Examples}

In this section we present two examples for the simulation of non-stationary stochastic processes by $3^{rd}$-order Spectral Representation Method. The first example involves the simulation of a stochastic process with separable time and frequency components which provides intuition to understand the POD based implementation. The second example considers the simulation of a fully non-stationary stochastic ground motion process using the Clough-Penzien power spectrum, which highlights the practical application of the simulation formula.

\subsection{Example 1: Non-stationary process with separable time and frequency contributions}

Let $f_s(t)$ represent a stationary stochastic process simulated using the $3^{rd}$-order Spectral Representation Method \cite{Vandanapu2021} from power spectrum $S(\omega)$ and bispectrum $B(\omega_1, \omega_2)$. Consider the evolutionary power spectral density takes the amplitude modulated form 
\begin{equation}
    S(t, \omega) = M^{2}(t) S(\omega)
\label{eqn:example_1_non_stationary_p}
\end{equation}
and the evolutionary bispectrum takes a similar amplitude modulated form given by
\begin{equation}
    B(t, \omega_1, \omega_2) = M(t)^{3} B(\omega_1, \omega_2).
\end{equation}
The resulting non-stationary process then takes the form
\begin{equation}
    f_{ns}(t) = M(t)f_{s}(t)
\end{equation}
The function $M(t)$ is simply a time-dependent modulating function which modifies the amplitude (but not the frequency content) of the stationary process to introduce non-stationarity.









\begin{figure}[!ht]
\centering
  \includegraphics[width=\linewidth]{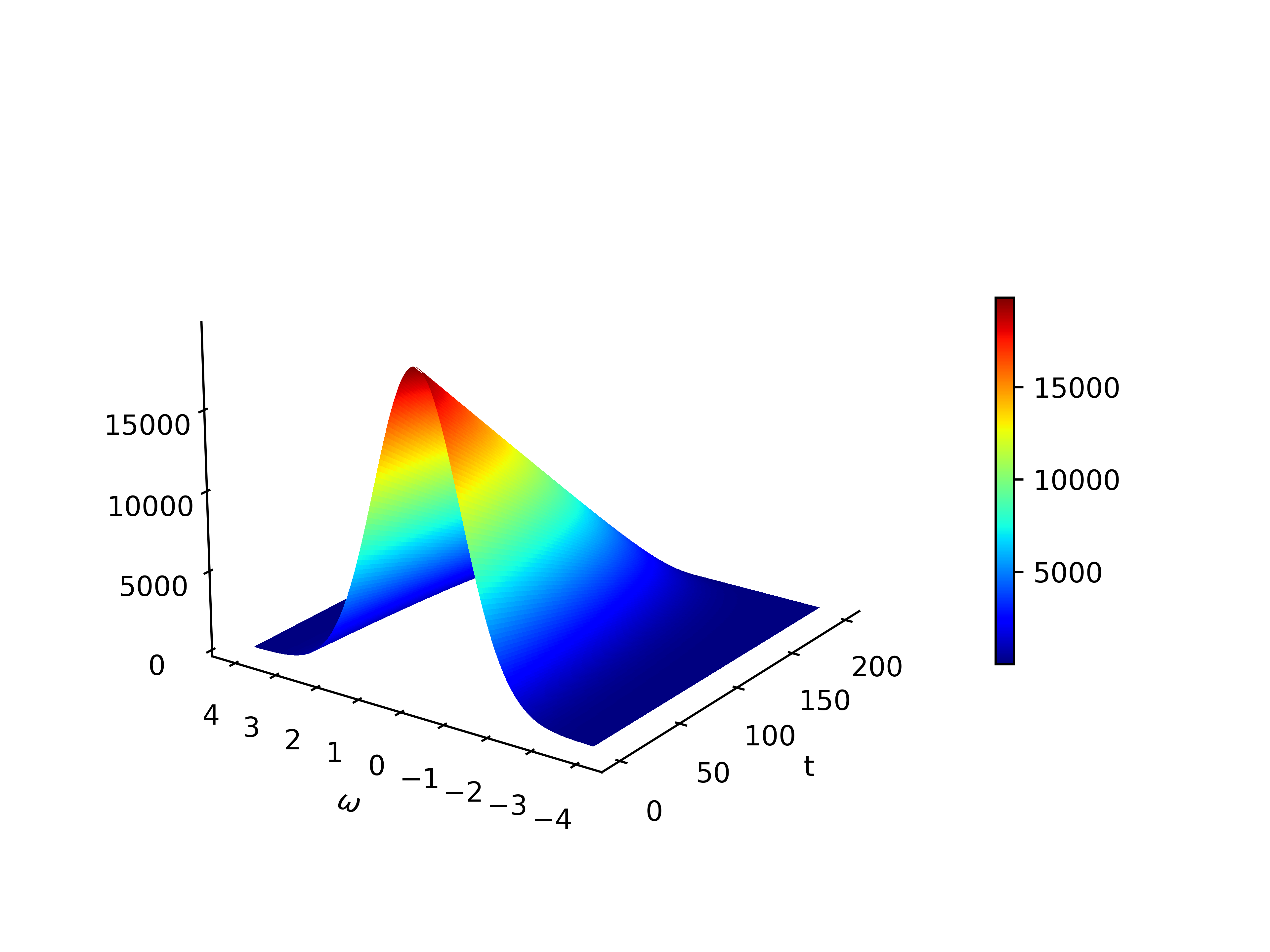}
  \caption{Example 1: Time-frequency separable evolutionary power spectral density}
  \label{fig:non_sta_pow_spec}
\end{figure}

Consider the specific evolutionary power spectrum of this form given by
\begin{equation}
    S(t, \omega) = 100(200 - t)e^{-\frac{1}{2}\omega^{2}}
\end{equation}
where $M(t) = \sqrt{200 - t}$ and $S(\omega) = 100e^{-\frac{1}{2}\omega^{2}}$ which is plotted in Figure \ref{fig:non_sta_pow_spec}.
The corresponding evolutionary bispectrum takes the form
\begin{equation}
    B(t, \omega_{1}, \omega_{2}) = \frac{2000}{3\sqrt{3(\omega_{1} + \omega_{2})}} (200 - t)^{\frac{3}{2}} e^{-\frac{1}{2}(\omega_{1}^{2} + \omega_{2}^{2} + \omega_{1}\omega_{2})}.
\end{equation}
Since $M(t) = \sqrt{200 - t}$, we can see that 
\begin{equation}
    B(\omega_1, \omega_2) = \frac{2000}{3\sqrt{3(\omega_{1} + \omega_{2})}} e^{-\frac{1}{2}(\omega_{1}^{2} + \omega_{2}^{2} + \omega_{1}\omega_{2})}
\end{equation}
The parameters used in the simulations are as follows
\begin{equation}
    T = 200 \text{sec} \quad  \omega_{0} = 4.02 \text{ rad/sec} \quad M = 256 \quad N = 128
\end{equation}

\begin{figure}[!ht]
\centering
\begin{subfigure}{.49\textwidth}
  \centering
  \includegraphics[width=0.99\linewidth]{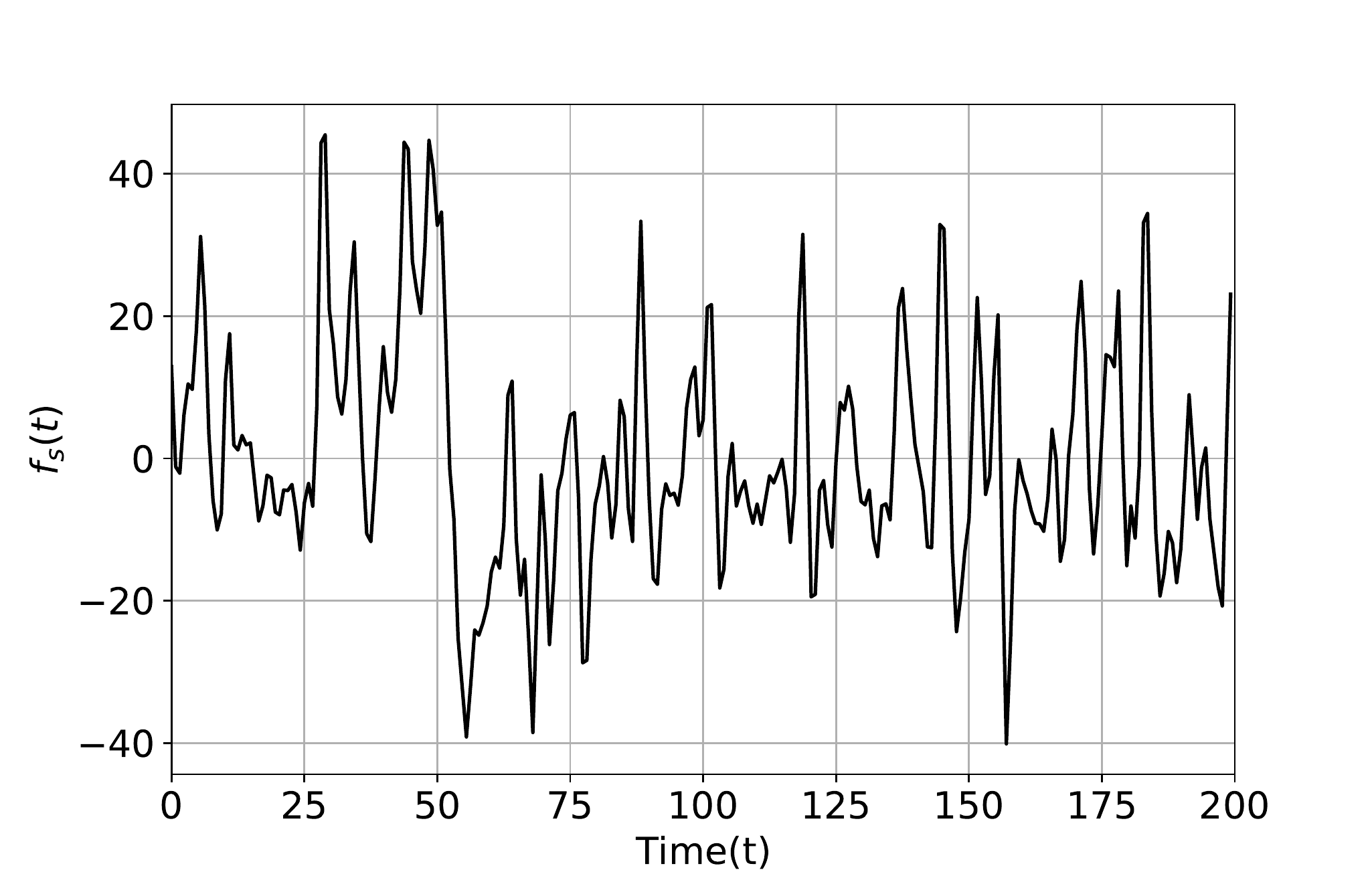}
  \caption{}
  \label{fig:example_1_sta_sample}
 \end{subfigure}
 \begin{subfigure}{.49\textwidth}
  \centering
  \includegraphics[width=0.99\linewidth]{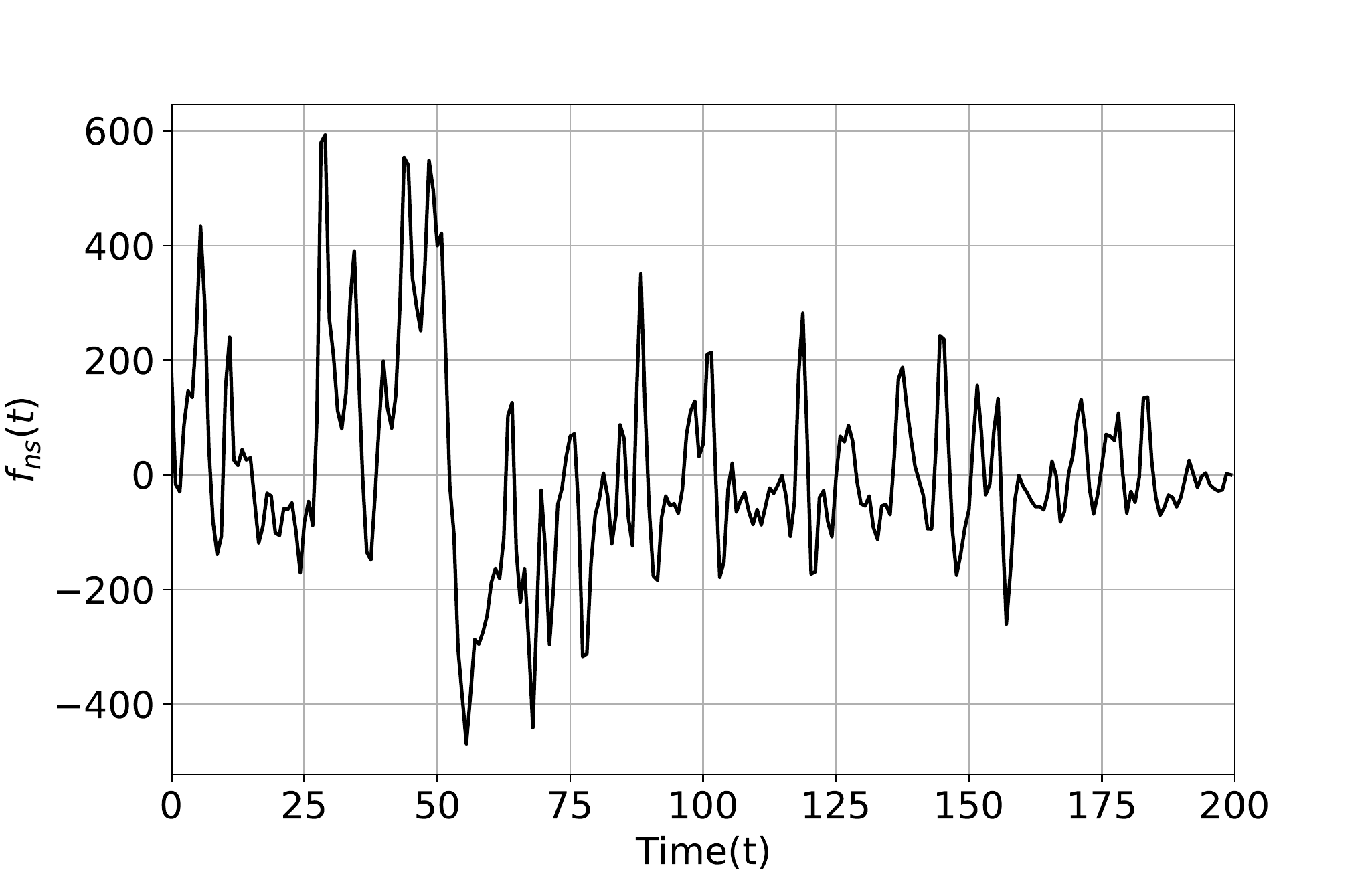}
  \caption{}
  \label{fig:example_1_non_sta_sample}
 \end{subfigure}
\caption{Example 1: Sample functions of (a) the stationary process $f_s(t)$ and (b) the amplitude modulated non-stationary processes $f_{ns}(t)$.}
\label{fig:example_1_samples}
\end{figure}

Sample functions of the stationary process, $f_{s}(t)$, and the non-stationary process, $f_{ns}(t)$, are plotted in Figure \ref{fig:example_1_samples}.

\begin{figure}[!ht]
    \centering
    \begin{subfigure}{.32\textwidth}
        \includegraphics[width=0.99\linewidth]{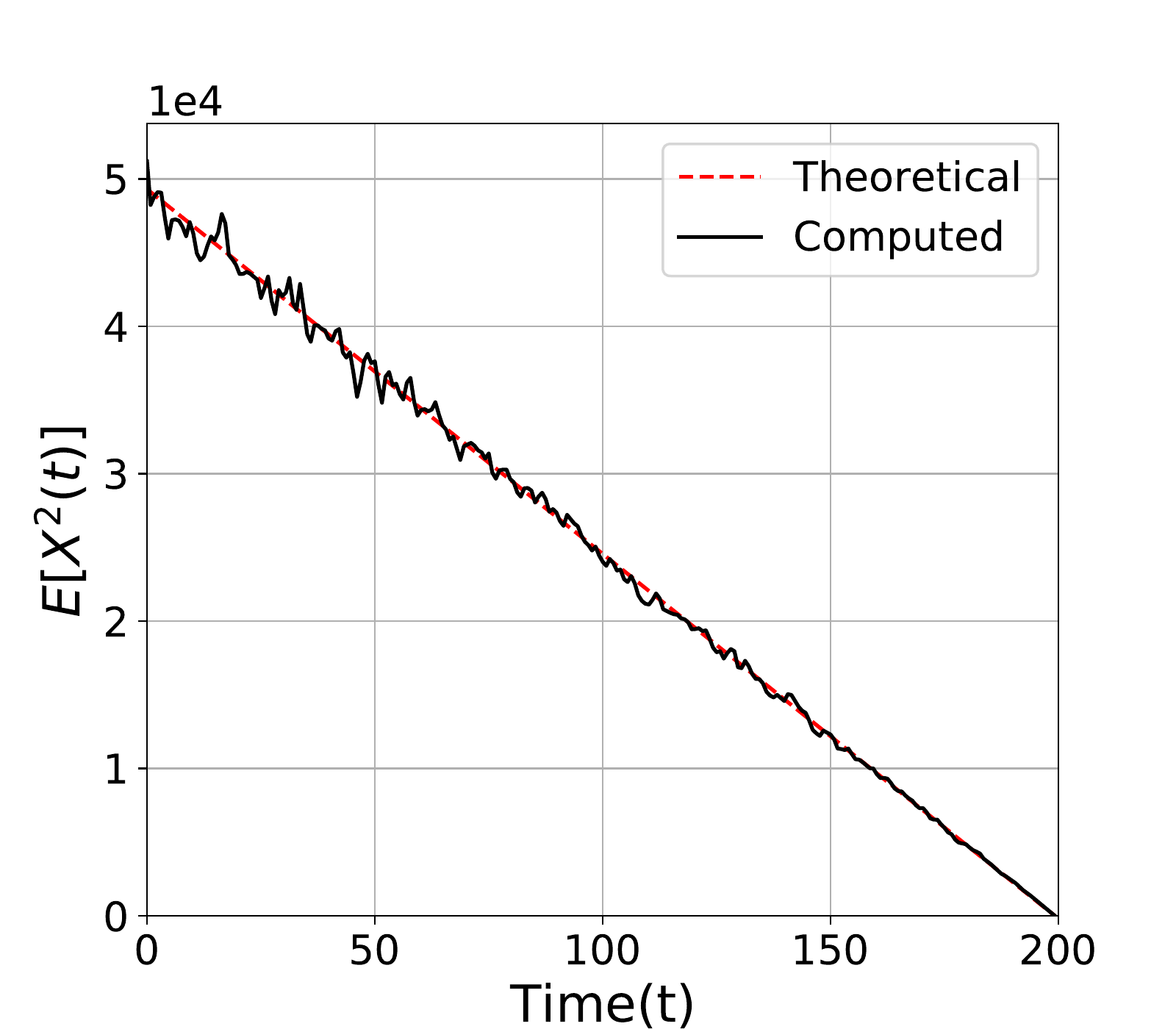}
        \caption{}
    \end{subfigure}
    \begin{subfigure}{.32\textwidth}
        \includegraphics[width=0.99\linewidth]{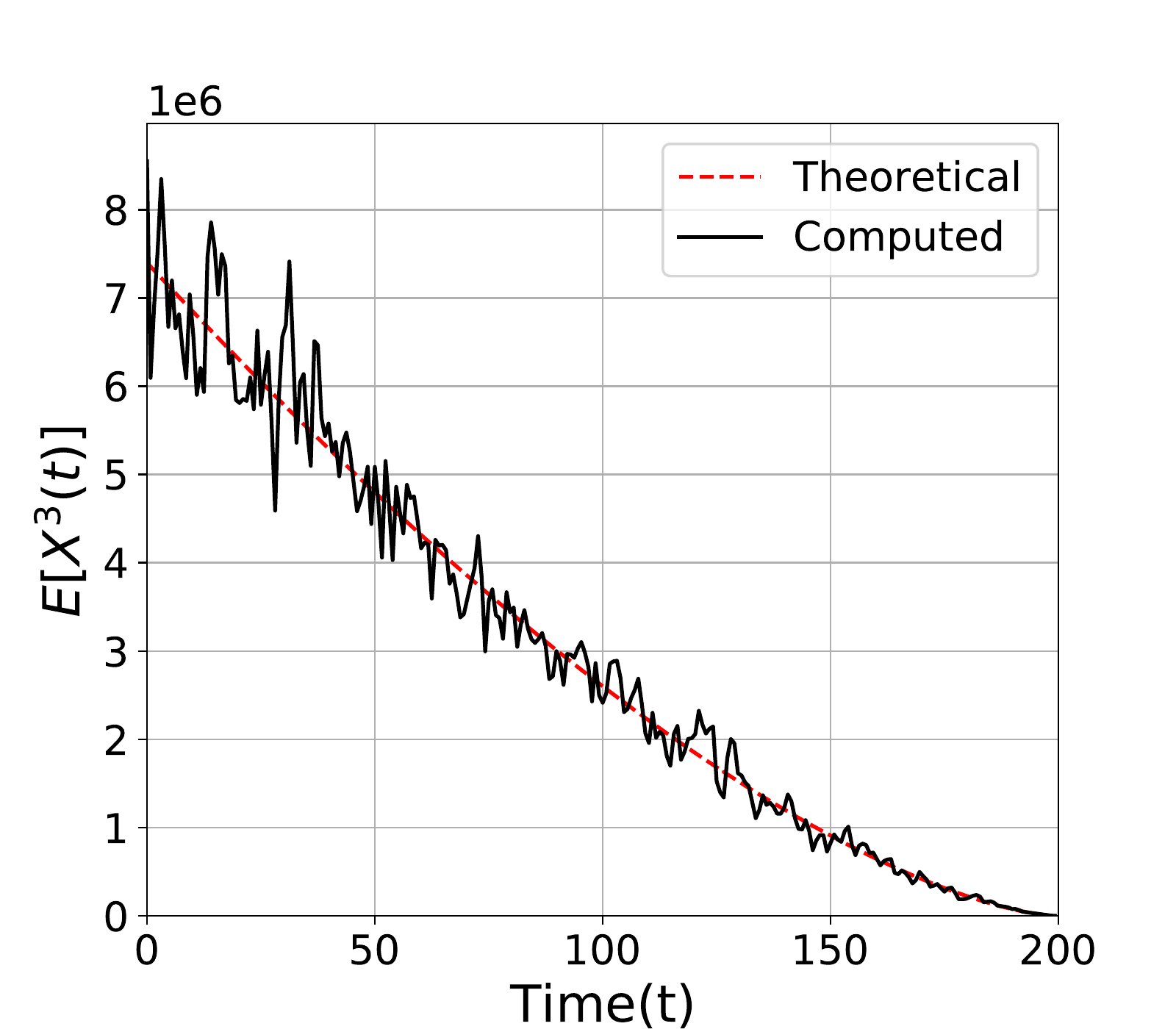}
        \caption{}
    \end{subfigure}
    \begin{subfigure}{.32\textwidth}
        \includegraphics[width=0.99\linewidth]{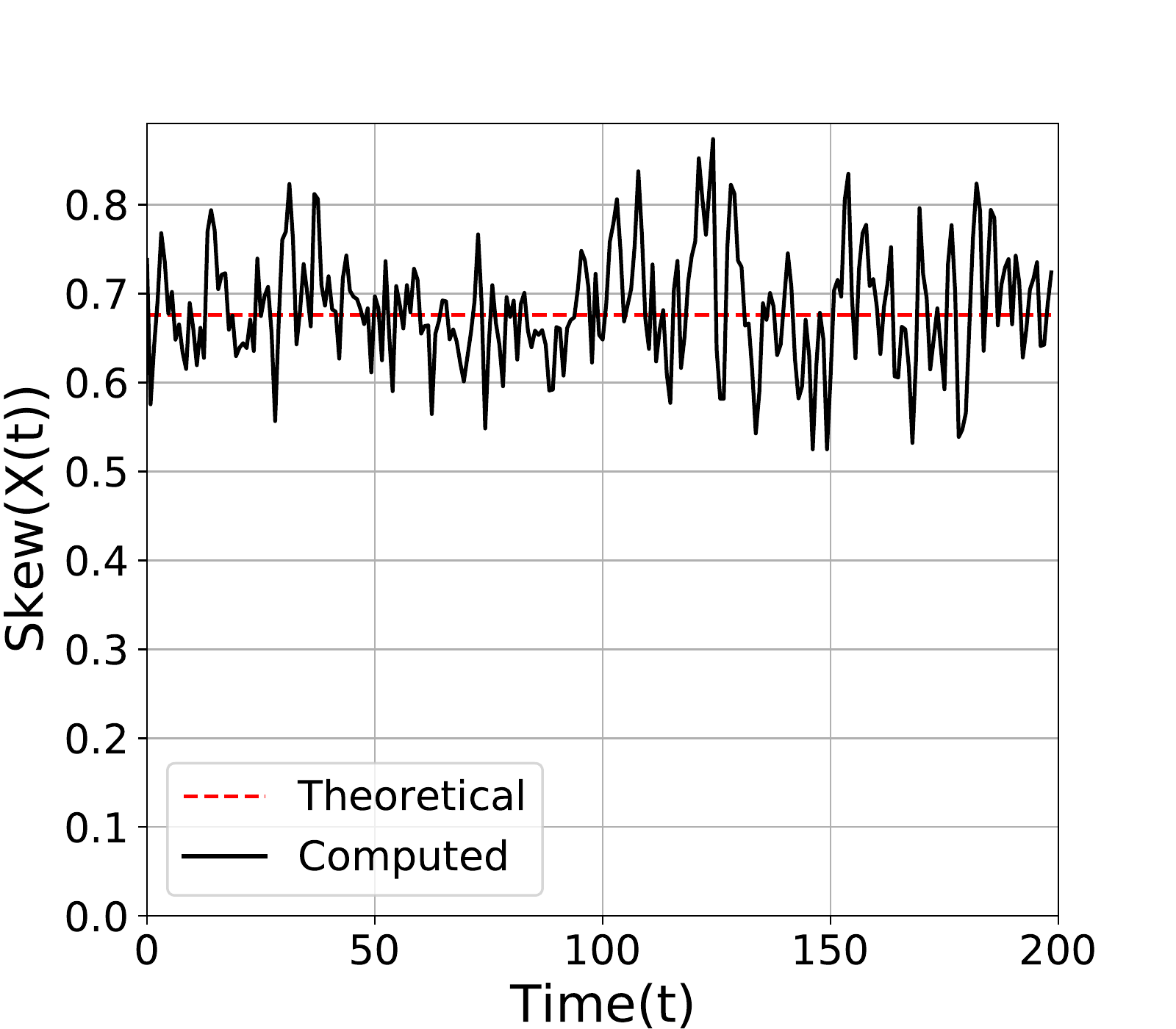}
        \caption{}
    \end{subfigure}
    \caption{Example 1: Time evolution of (a) the second moment, (b) the third moment, and (c) the skewness for the separable non-stationary stochastic process.}
    \label{fig:example_1_statistics_time}
\end{figure}

Figure \ref{fig:example_1_statistics_time} shows the time-evolution of the second moment, third moment, and skewness of the process computed from 10,000 sample non-stationary processes along with theoretical values derived from the evolutionary spectra. We can see that the properties of the simulated non-stationary processes closely match the theoretical values.

\begin{figure}[!ht]
\centering
  \includegraphics[width=0.5\linewidth]{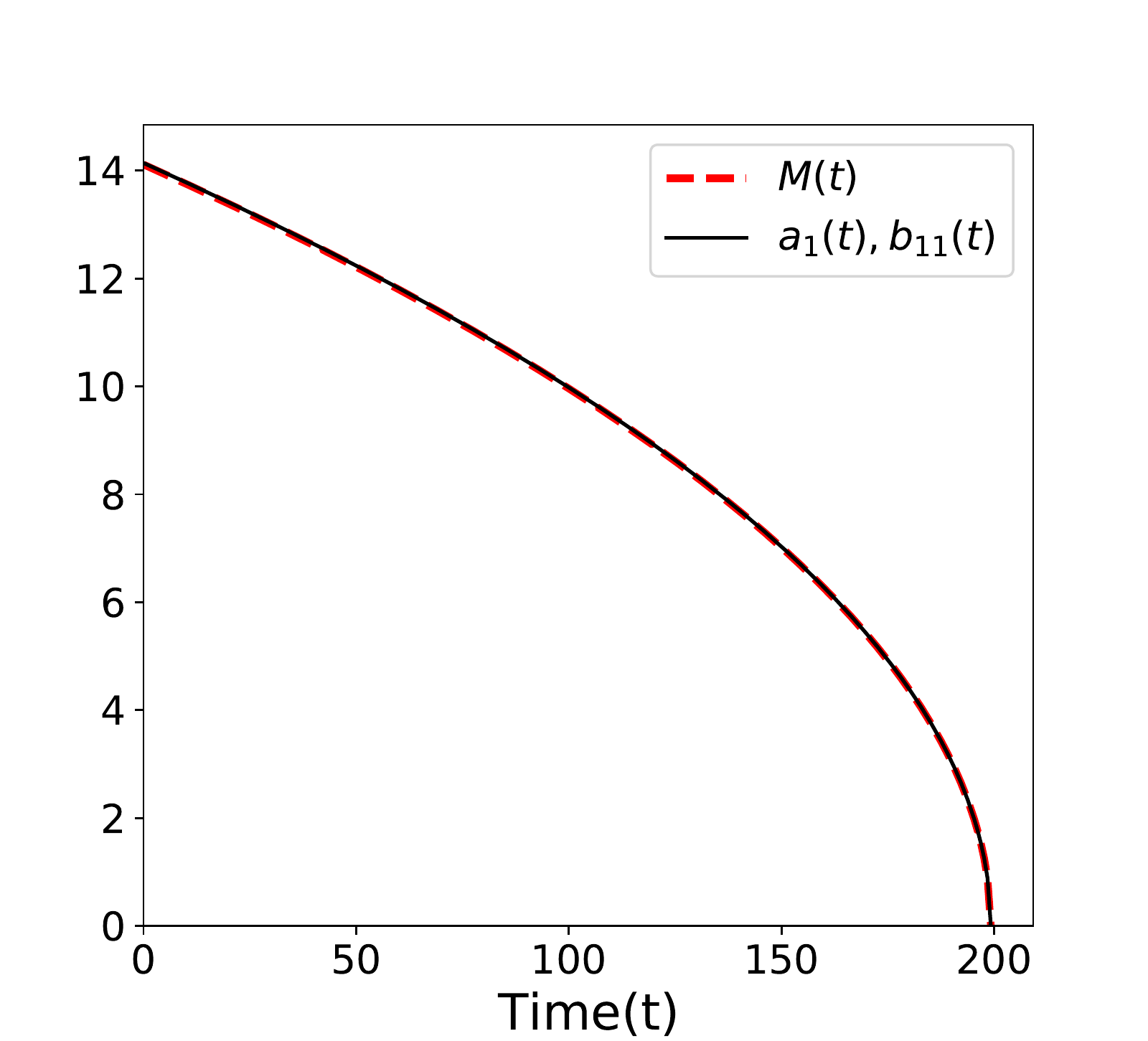}
  \caption{Example 1: Modulating Function $M(t)$ and the first time-dependent amplitude functions $a_{1}(t), b_{11}(t)$ from the POD for the separable evolutionary power spectrum.}
  \label{fig:example_1_modulating_function}
\end{figure}

Proper Orthogonal Decomposition of the evolutionary power spectrum reveals that the first component in the pure component of the expansion $a_{1}(t)$ and the interactive component of the expansion $b_{11}(t)$ are both equal to the modulating function $M(t)$ as seen in Figure \ref{fig:example_1_modulating_function}. This demonstrates that the POD can successfully decompose a separable evolutionary power spectra with only a single orthogonal component. This further illustrates the interpretation that the POD based implementation is akin to finding a series of modulating functions and underlying stationary spectra spectra.


\subsection{Non-stationary seismic ground motion}

\begin{figure}[!ht]
\centering
  \includegraphics[width=\linewidth]{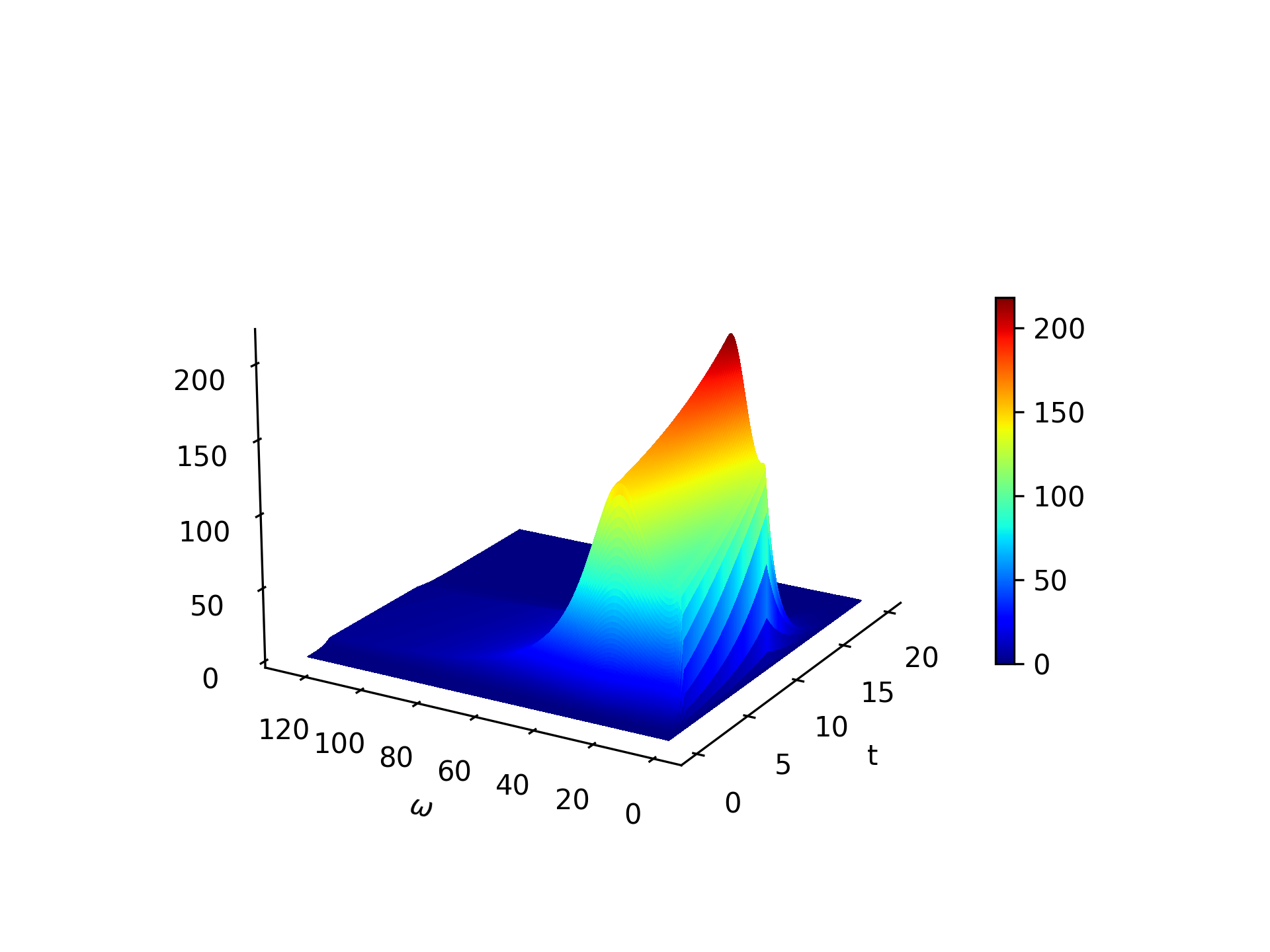}
  \caption{Example 2: Clough-Penzien Evolutionary Spectrum}
  \label{fig:clough_penzien_power_spectrum}
\end{figure}

Next, we consider the simulation of a non-stationary stochastic ground motion process using the Kanai-Tajimi evolutionary power spectrum with Clough-Penzien correction \cite{clough_dynamics_1975}, which contains both frequency and amplitude modulation. The equation for the Kanai-Tajimi evolutionary spectrum is given by
\begin{equation}
    S_{kt}(t, \omega) = \frac{1 + 4\zeta_g^{2}\big(\frac{\omega}{\omega_{g}}\Big)^{2}}{\Big[1 - \Big(\frac{\omega}{\omega_g}\Big)^{2} \Big]^{2} + 4\zeta_g^{2}\big(\frac{\omega}{\omega_{g}}\Big)^{2}},
\end{equation}
and the Clough-Penzien correction factor is defined as
\begin{equation}
    \gamma_{cp}(t, \omega) = \frac{\big(\frac{\omega}{\omega_{f}}\Big)^{4}}{\Big[1 - \Big(\frac{\omega}{\omega_f}\Big)^{2} \Big]^{2} + 4\zeta_f^{2}\big(\frac{\omega}{\omega_{f}}\Big)^{2}}
\end{equation}
where the time-dependent, non-stationary parameters are defined as
\begin{equation}
\begin{aligned}
    & \omega_g = 30 - 1.25t \\
    & \omega_f = 0.1 \omega_g \\
    & \zeta_g = 0.5 + 0.005t \\
    & \zeta_f = 0.1 \zeta_g \\
\end{aligned}
\end{equation}
such that $\omega_g$, $\zeta_g$ are the characteristic frequency and damping of the ground and $\omega_f$, $\zeta_f$ are the filtering parameters of the Clough-Penzien correction. The resulting Clough-Penzien evolutionary power spectrum is given by
\begin{equation}
    S_{cp}(t, \omega) = S_{kt}(t, \omega)\gamma_{cp}(t, \omega)
\label{eqn:clough_penzien_ps}
\end{equation}
and is plotted in Figure \ref{fig:clough_penzien_power_spectrum}.
We define the associated bispectrum of the non-stationary process as
\begin{equation}
    B(t, \omega_i, \omega_j) = \frac{2\sqrt{S_{cp}(t, \omega_i)S_{cp}(t, \omega_j)S_{cp}(t, \omega_i + \omega_j)}}{3\sqrt{3(\omega_i + \omega_j)}}
\end{equation}

\begin{figure}[!ht]
    \centering
    \begin{subfigure}{.45\textwidth}
        \centering
        \includegraphics[width=\linewidth]{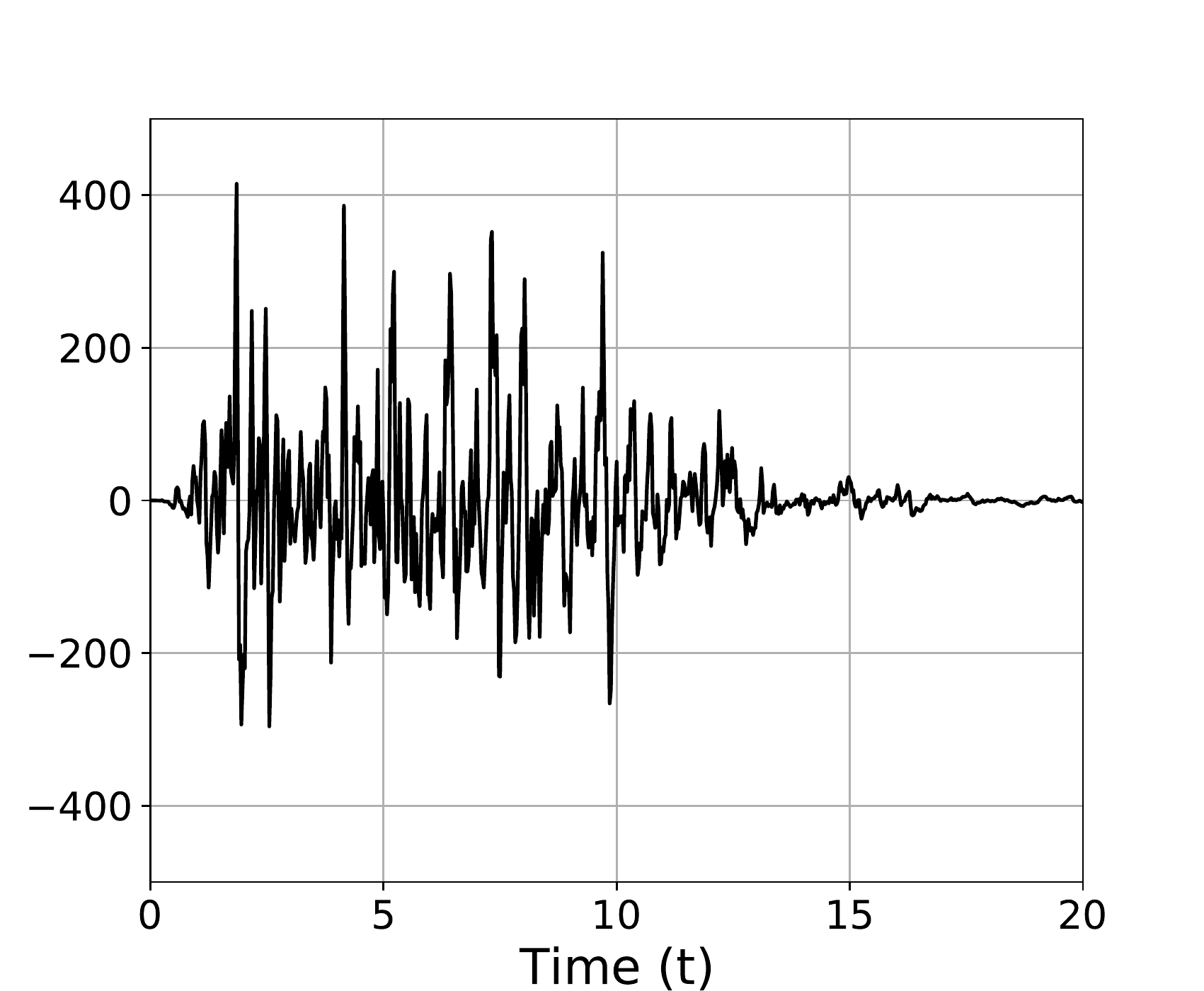}
        \caption{}
        \label{fig:sample_bsrm}
    \end{subfigure}
        \begin{subfigure}{.45\textwidth}
        \centering
        \includegraphics[width=\linewidth]{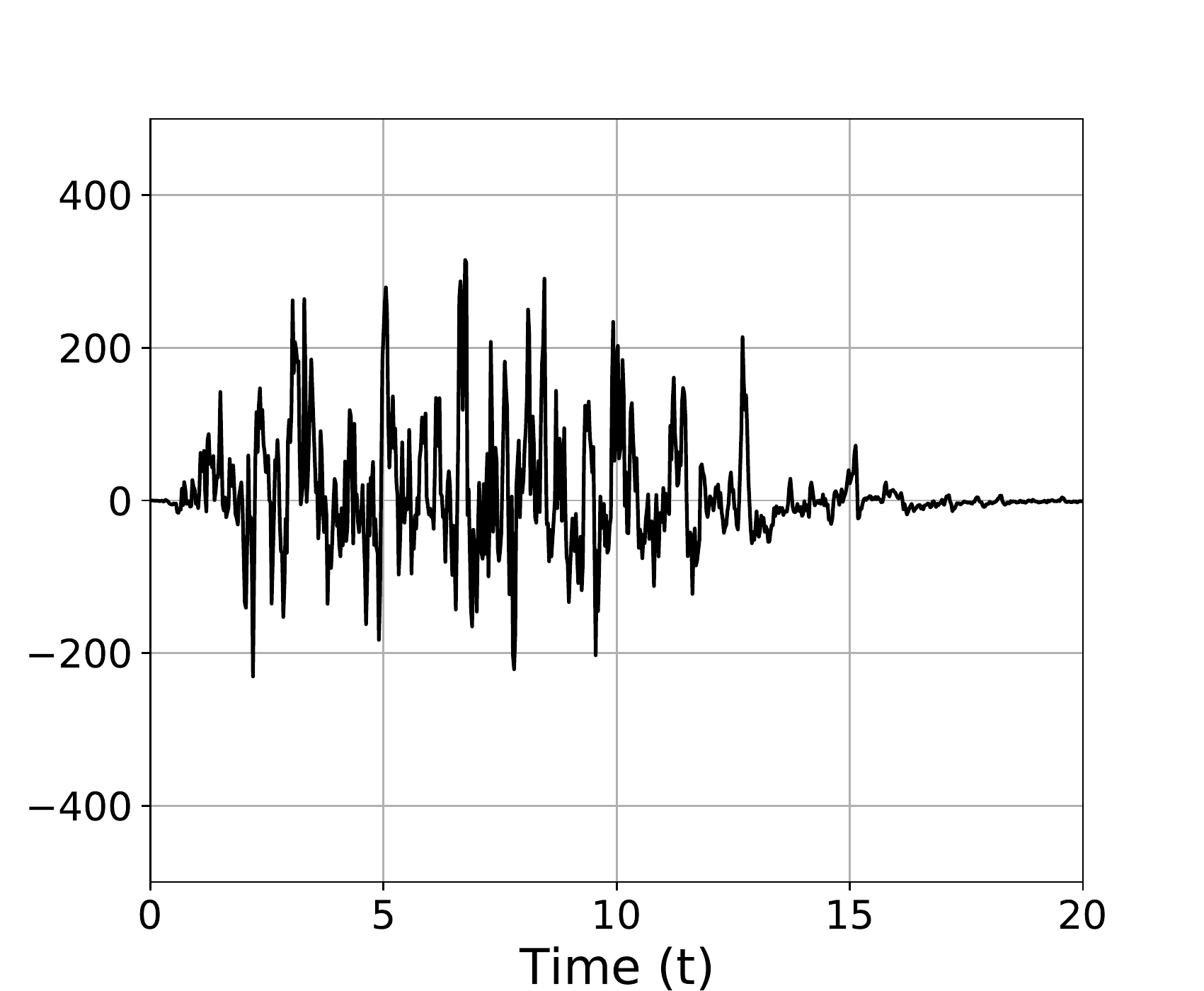}
        \caption{}
        \label{fig:sample_bsrm_pod}
    \end{subfigure}
\caption{Example 2: Representative sample functions of the third-order non-stationary stochastic ground motion process simulated using (a) the direct sum of cosines and (b) the POD implementation.}
\label{fig:samples_bsrm}
\end{figure}

The non-stationary process is simulated using both the direct sum of cosines and the POD-based implementation. The parameters of the simulation are given by
\begin{equation}
    T = 20 \text{}  \quad  \omega_{0} = 125.66 \text{ rad/sec} \quad M = 800 \quad N = 400 \quad n_q = 10
\end{equation}
Sample function plots for each implementation are shown in Figure \ref{fig:samples_bsrm}, demonstrating that both methods are capable of simulating the random process.

\begin{figure}[!ht]
    \centering
    \begin{subfigure}{.32\textwidth}
        \centering
        \includegraphics[width=\linewidth]{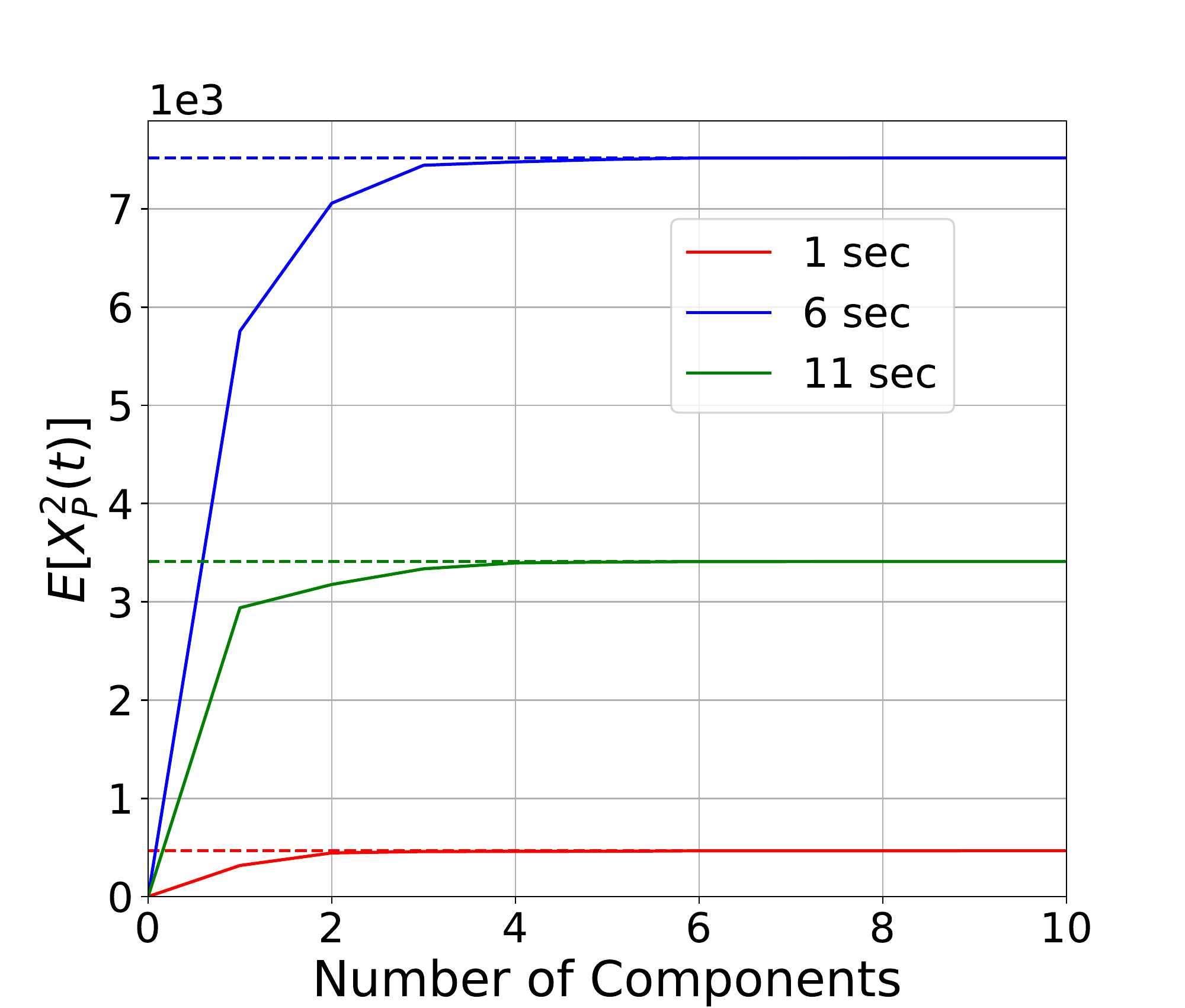}
        \caption{}
        \label{fig:pure_variance_convergence}
    \end{subfigure}
    \begin{subfigure}{.32\textwidth}
        \centering
        \includegraphics[width=\linewidth]{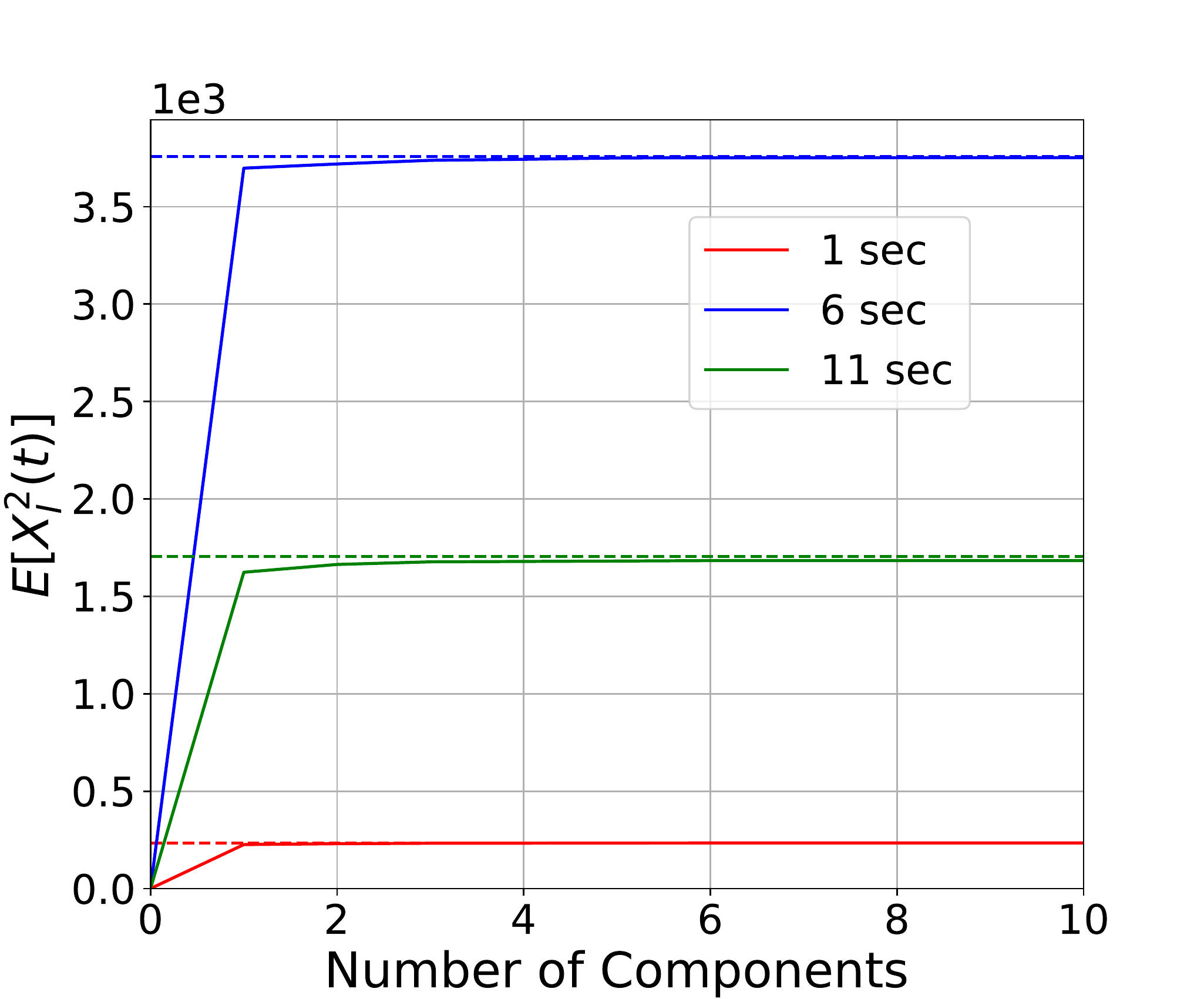}
        \caption{}
        \label{fig:interactive_variance_convergence}
    \end{subfigure}
    \begin{subfigure}{.32\textwidth}
        \centering
        \includegraphics[width=\linewidth]{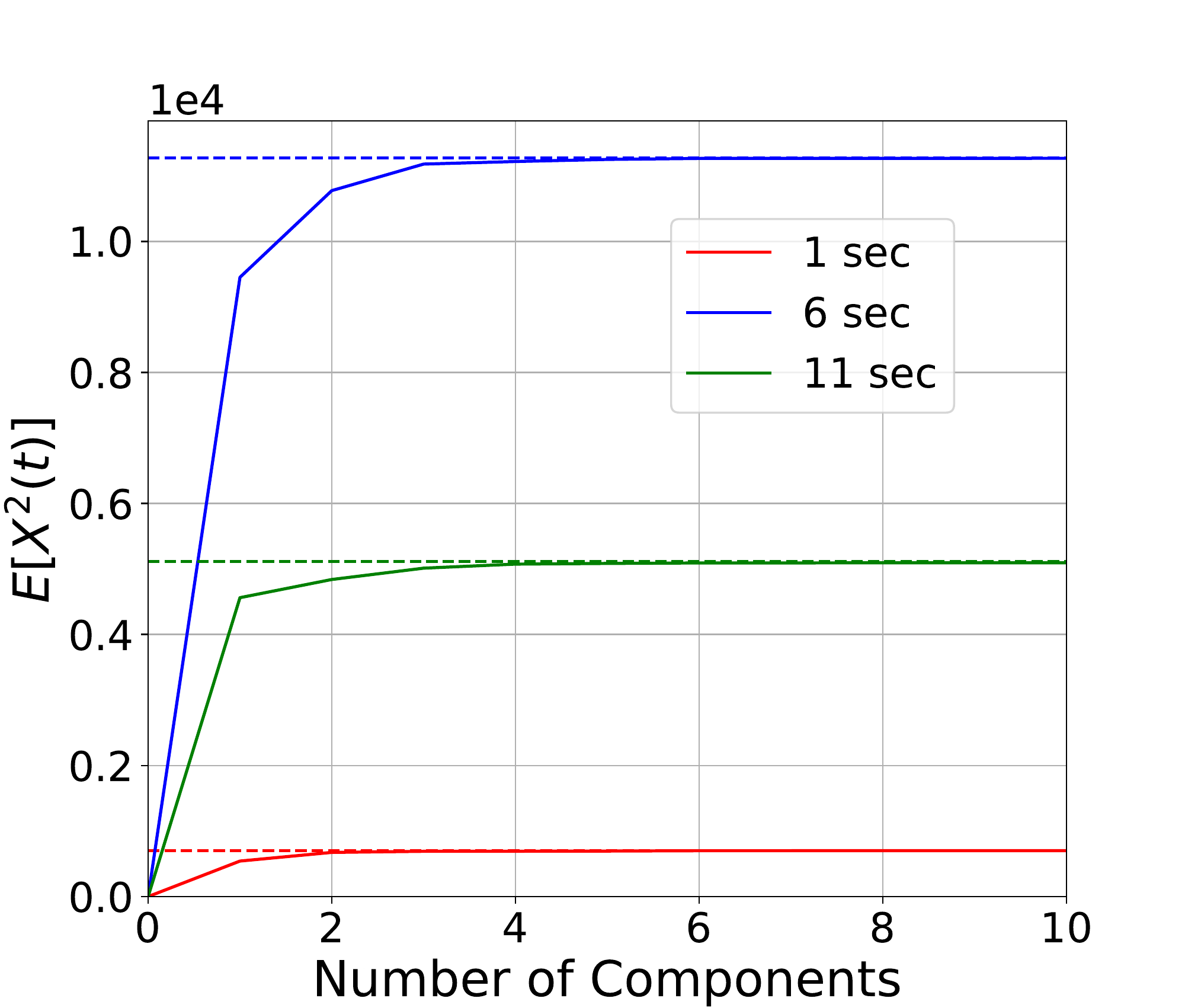}
        \caption{}
        \label{fig:variance_convergence}
    \end{subfigure}
\caption{Example 2: Convergence of (a) the pure component of the second moment, (b) the interactive component of the second moment, and (c) the total second moment at different time instances with increasing number of POD components.}
\label{fig:second_order_statistics_convergence}
\end{figure}

\begin{figure}[!ht]
    \centering
    \begin{subfigure}{.49\textwidth}
        \centering
        \includegraphics[width=\linewidth]{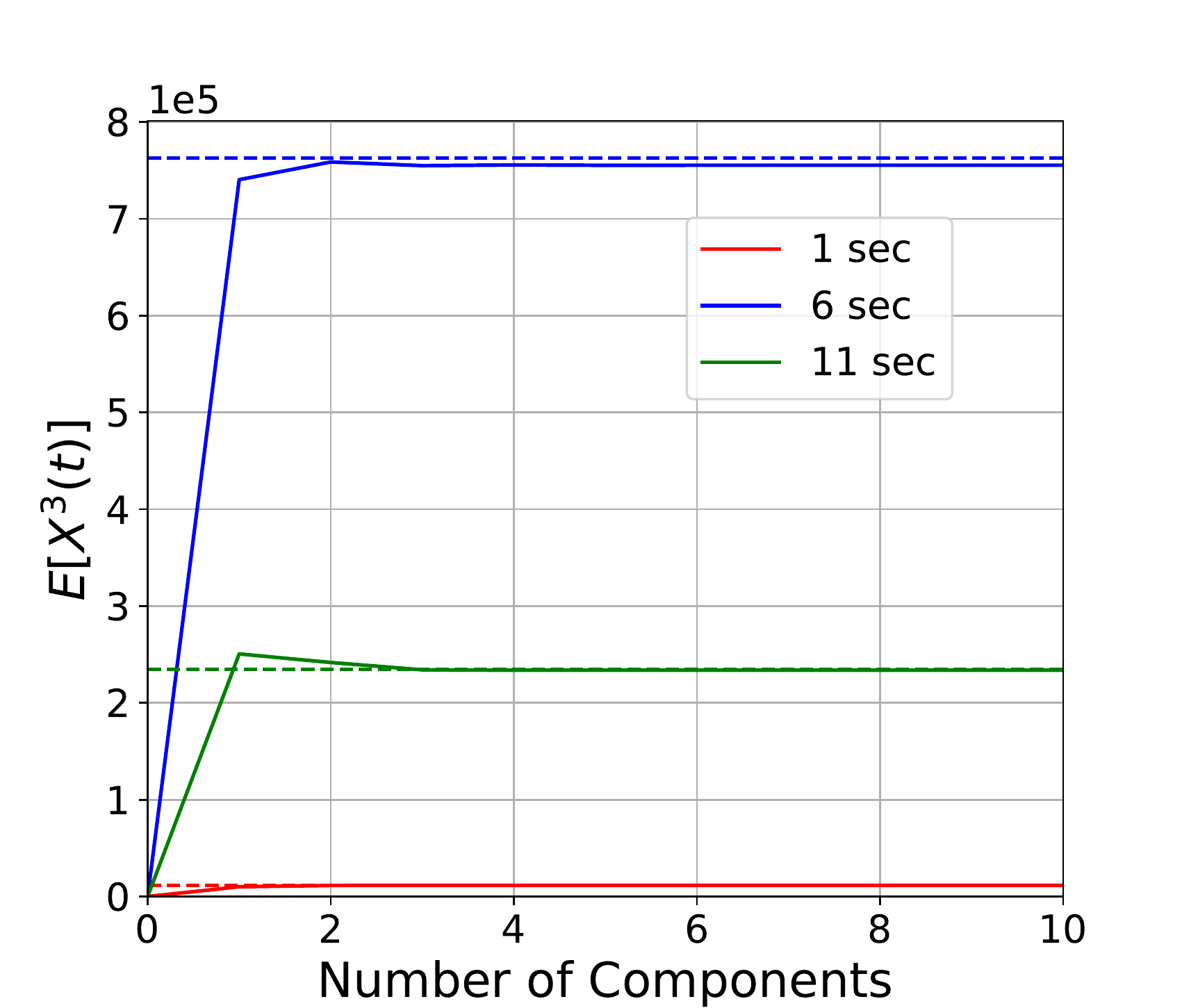}
        \caption{}
        \label{fig:third_order_convergence}
    \end{subfigure}
    \begin{subfigure}{.49\textwidth}
        \centering
        \includegraphics[width=\linewidth]{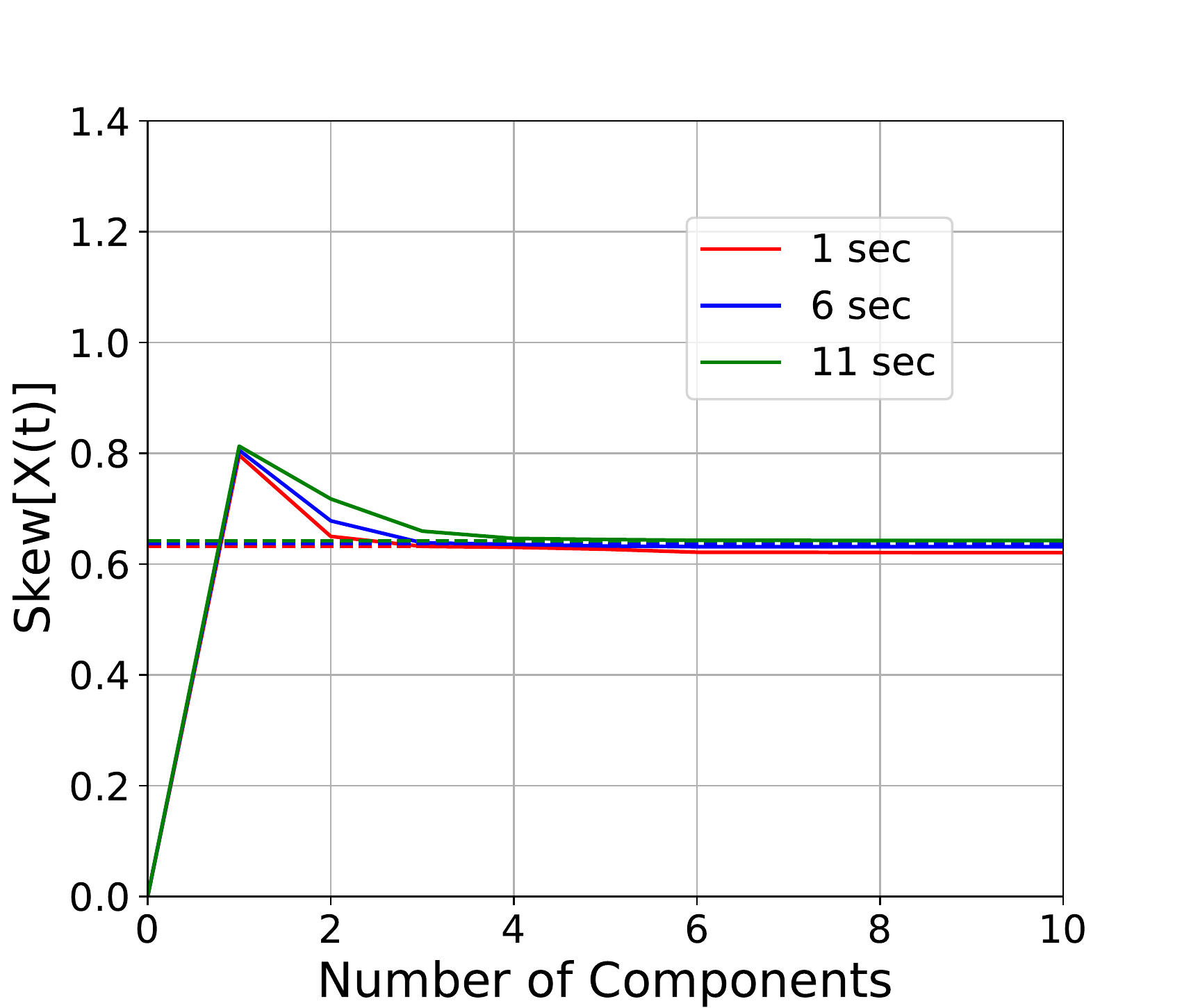}
        \caption{}
        \label{fig:skewness_convergence}
    \end{subfigure}
\caption{Example 2: Convergence of (a) the third moment, and (b) the skewness at different time instances with increasing number of POD components.}
\label{fig:third_order_statistics_convergence}
\end{figure}

\begin{figure}[!ht]
    \centering
    \begin{subfigure}{.32\textwidth}
        \includegraphics[width=0.99\linewidth]{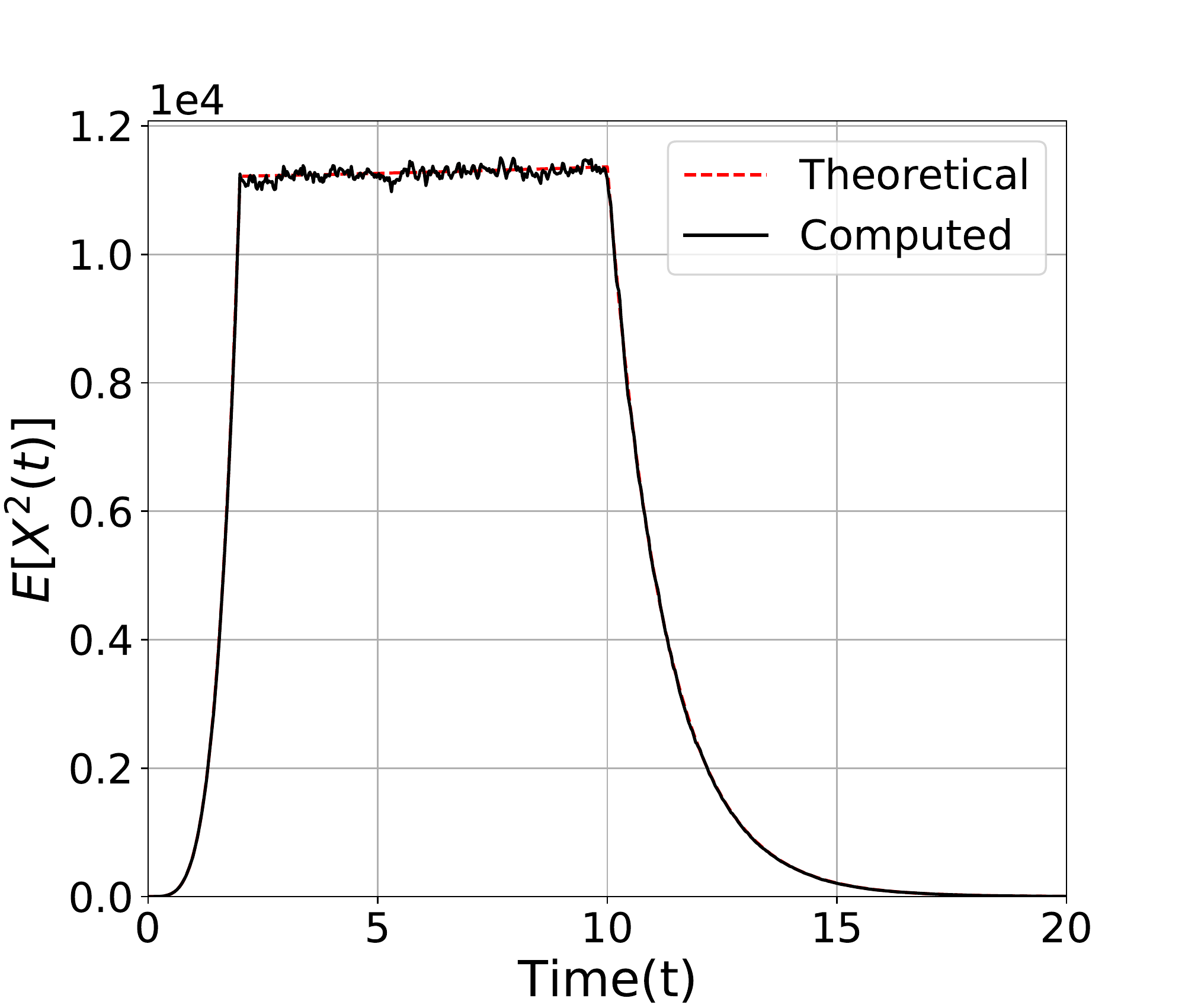}
        \caption{}
    \end{subfigure}
    \begin{subfigure}{.32\textwidth}
        \includegraphics[width=0.99\linewidth]{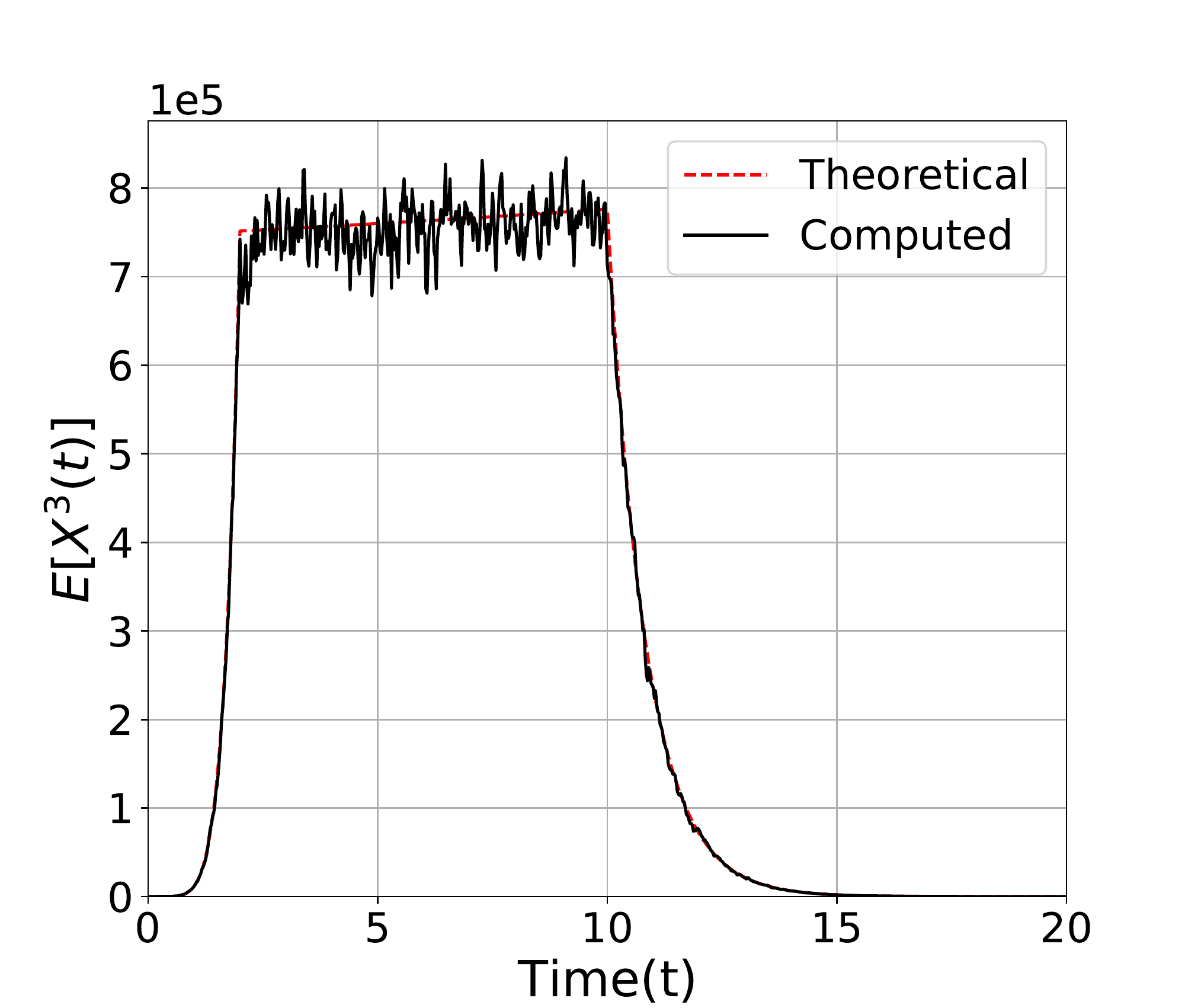}
        \caption{}
    \end{subfigure}
    \begin{subfigure}{.32\textwidth}
        \includegraphics[width=0.99\linewidth]{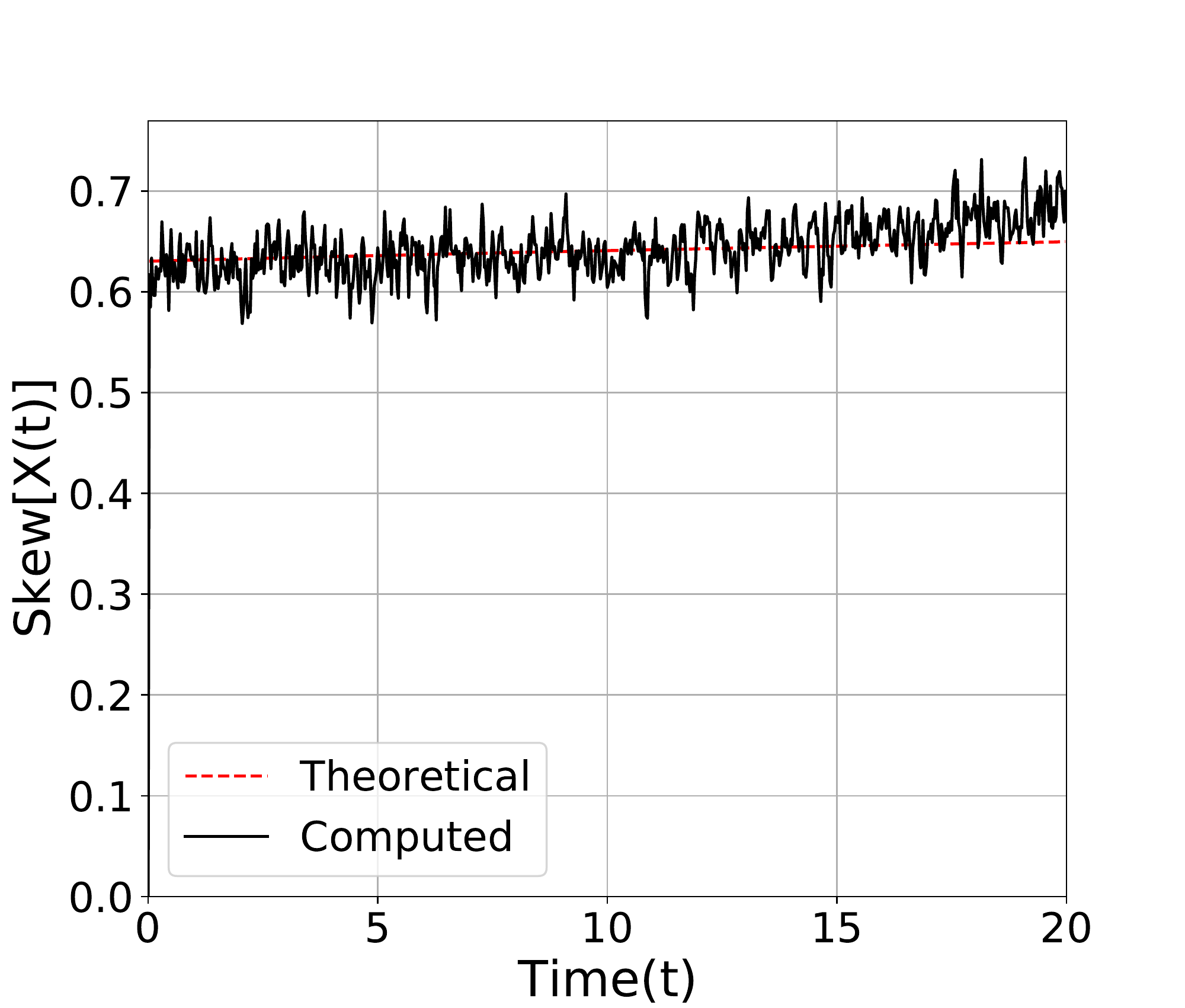}
        \caption{}
    \end{subfigure}
    \caption{Example 2: Evolution of second-order moment(variance), third-order moment and skewness with time.}
    \label{fig:example_2_statistics_time}
\end{figure}

Sample statistics of the simulated process at three different time instants from 10,000 simulations using both the sum of cosines and POD implementation with increasing number of POD components are shown in Figures \ref{fig:second_order_statistics_convergence}, \ref{fig:third_order_statistics_convergence}. The dotted lines in these plots represents the theoretical values, while the solid lines represent the simulation results.
We can see that the statistics converge very rapidly with only 4 POD components. We further plot the time-varying second and third moments of the simulated process, estimated from 100,000 samples, using the first 4 POD modes in Figure \ref{fig:example_2_statistics_time} along with their theoretical values.
Here, we can see that the statistical properties of the process match the their theoretical moments with high accuracy throughout the time duration of the process.

Next, we compare the computation time for simulation of non-stationary processes by the $2^{nd}$-order and $3^{rd}$-order Spectral Representation Method using both the standard sum of cosines approach and the POD methods for increasing number of samples. Results are shown in Tables \ref{table:2_order_nsampels} and \ref{table:3_order_nsampels} for the  $2^{nd}$-order and $3^{rd}$-order SRM, respectively.
\begin{table}[!ht]
    \centering
    \begin{tabular}{|c|c|c|c|c|}
        \hline
        & \multicolumn{4}{c|}{\textbf{Time (sec)}} \\
        \hline
        \textbf{Number of Samples} & \textbf{Sum of Cosines} & \multicolumn{3}{c|}{\textbf{POD}} \\ \cline{3-5}
        & & Total & Decomposition & Simulation \\
        \hline
        1 & 3.0971 & 0.7885 & 0.7864 & 0.0021 \\
        10 & 3.4194 & 0.7705 & 0.7614 & 0.0091 \\
        100 & 5.1791 & 0.8472 & 0.7590 & 0.0882 \\
        1000 & 29.0474 & 1.7306 & 0.7603 & 0.9703 \\
        10000 & 250.0017 & 16.5228 & 0.7630 & 15.7598 \\
        100000 & 2738.2145 & 137.4785 & 0.7633 & 136.7152 \\
        \hline
    \end{tabular}
    \caption{Comparison of computation time for the sum of cosines formula and the POD method for simulating $2^{nd}$-order non-stationary stochastic ground motion processes.}
    \label{table:2_order_nsampels}
\end{table}
\begin{table}[!ht]
    \centering
    \begin{tabular}{|c|c|c|c|c|}
        \hline
        & \multicolumn{4}{c|}{\textbf{Time (sec)}} \\
        \hline
        \textbf{Number of Samples} & \textbf{Sum of Cosines} & \multicolumn{3}{c|}{\textbf{POD}} \\ \cline{3-5}
        & & Total & Decomposition & Simulation \\
        \hline
        1 & 435.1617 & 36.7289 & 22.9394 & 13.7895 \\
        10 & 480.0146 & 39.0056 & 23.4440 & 15.5616 \\
        100 & 659.3620 & 39.3217 & 23.3529 & 15.9688	\\
        1000 & 3322.8135 & 71.7228 & 23.5590 & 48.1638 \\
        \hline
    \end{tabular}
    \caption{Comparison of the computation time for the sum of cosines formula and the POD method for simulating $3^{rd}$-order non-stationary stochastic ground motion processes.}
    \label{table:3_order_nsampels}
\end{table}
The tables show total CPU time for both the sum of cosines approach and the POD based implementation, where times for POD based implementation are broken down into three components: total CPU time, CPU time for the POD of the evolutionary power spectrum, and CPU time for simulation of the stationary samples using FFT. We can see that the time for the POD remains constant, while the FFT simulation scales with the number of simulations. Hence, the computation time will scale with sample size and with $N_q$ (i.e. the number of FFTs required for a single simulation). More specifically, simulation time is proportional to the square of $N_q$ in the POD method (see  Eq.\ \eqref{eqn:pod-time-complexity}), so choosing fewer components helps drastically reduce the simulation time.The tables further show that for the given $M$ and $N$, i.e. the given number of time and frequency intervals, using the POD based implementation is always computationally beneficial for both $2^{nd}$ and $3^{rd}$-order processes, even when simulating only a single sample function. Note that, in Table \ref{table:3_order_nsampels}, we do not simulate more than 1000 realizations because the sum of cosines implementation become intractable.


\section{Conclusion}

In this paper, the $3^{rd}$-order Spectral Representation Method has been extended for the simulation of non-stationary stochastic processes. First, the standard sum of cosines form was derived. This conventional implementation is computationally prohibitive. To alleviate the computational burden, a Proper Orthogonal Decomposition (POD) based implementation was presented, which enables the use of Fast Fourier Transform to significant speed up the simulations. Two example processes are simulated to highlight the advantages of the proposed method. The first example considers a time-frequency separable process, where we show that the POD implementation correctly separates the components for simulation. The second example considers an inseparable ground motion process, demonstrating convergence of the POD implementation in the statistical response using only a small number of components. The second example further explores the computational gains afforded by the POD implementation and demonstrates that the POD implementation can be orders of magnitude more efficient and enable the simulation of processes that are infeasible to simulate with the standard sum of cosines implementation.

\section{Acknowledgement}

This work has been supported by the National Science Foundation under award number 1652044.

\bibliographystyle{model1-num-names}
\bibliography{elsarticle-template-1-num.bib}

\begin{thebibliography}{32}
\expandafter\ifx\csname natexlab\endcsname\relax\def\natexlab#1{#1}\fi
\providecommand{\bibinfo}[2]{#2}
\ifx\xfnm\relax \def\xfnm[#1]{\unskip,\space#1}\fi
\bibitem[{Huang et~al.(2001)Huang, Quek, and Phoon}]{Huang2001}
\bibinfo{author}{S.~P. Huang}, \bibinfo{author}{S.~T. Quek},
  \bibinfo{author}{K.~K. Phoon},
\newblock \bibinfo{title}{{Convergence study of the truncated Karhunen–Loeve
  expansion for simulation of stochastic processes}},
\newblock \bibinfo{journal}{International Journal for Numerical Methods in
  Engineering} \bibinfo{volume}{52} (\bibinfo{year}{2001})
  \bibinfo{pages}{1029--1043}.
\bibitem[{Shinozuka(1972)}]{Shinozuka1972}
\bibinfo{author}{M.~Shinozuka},
\newblock \bibinfo{title}{{Monte Carlo Solution of Structural Dynamics}},
\newblock \bibinfo{journal}{Computers and Structures} \bibinfo{volume}{2}
  (\bibinfo{year}{1972}) \bibinfo{pages}{855--874}.
\bibitem[{Yang(1972)}]{Yang1972}
\bibinfo{author}{J.-N. Yang},
\newblock \bibinfo{title}{{Simulation of random envelope processes}},
\newblock \bibinfo{journal}{Journal of Sound and Vibration}
  \bibinfo{volume}{21} (\bibinfo{year}{1972}) \bibinfo{pages}{73--85}.
\bibitem[{Deodatis(1996)}]{Deodatis1996a}
\bibinfo{author}{G.~Deodatis},
\newblock \bibinfo{title}{{Simulation of Ergodic Multivariate Stochastic
  Processes}},
\newblock \bibinfo{journal}{Journal of Engineering Mechanics}
  \bibinfo{volume}{122} (\bibinfo{year}{1996}) \bibinfo{pages}{778--787}.
\bibitem[{Deodatis and Shinozuka(1989)}]{Deodatis1989}
\bibinfo{author}{G.~Deodatis}, \bibinfo{author}{M.~Shinozuka},
\newblock \bibinfo{title}{{Simulation of Seismic Ground Motion Using Stochastic
  Waves}},
\newblock \bibinfo{journal}{Journal of Engineering Mechanics}
  \bibinfo{volume}{115} (\bibinfo{year}{1989}) \bibinfo{pages}{2723--2737}.
\bibitem[{Shinozuka and Deodatis(1996)}]{Shinozuka1996}
\bibinfo{author}{M.~Shinozuka}, \bibinfo{author}{G.~Deodatis},
\newblock \bibinfo{title}{{Simulation of Multi-Dimensional Gaussian Stochastic
  Fields by Spectral Representation}},
\newblock \bibinfo{journal}{Applied Mechanics Reviews} \bibinfo{volume}{49}
  (\bibinfo{year}{1996}) \bibinfo{pages}{29--53}.
\bibitem[{Grigoriu(1993)}]{grigoriu1993spectral}
\bibinfo{author}{M.~Grigoriu},
\newblock \bibinfo{title}{On the spectral representation method in simulation},
\newblock \bibinfo{journal}{Probabilistic Engineering Mechanics}
  \bibinfo{volume}{8} (\bibinfo{year}{1993}) \bibinfo{pages}{75--90}.
\bibitem[{Puig et~al.(2002)Puig, Poirion, and Soize}]{Puig2002}
\bibinfo{author}{B.~Puig}, \bibinfo{author}{F.~Poirion},
  \bibinfo{author}{C.~Soize},
\newblock \bibinfo{title}{{Non-Gaussian simulation using Hermite polynomial
  expansion: Convergences and algorithms}},
\newblock \bibinfo{journal}{Probabilistic Engineering Mechanics}
  \bibinfo{volume}{17} (\bibinfo{year}{2002}) \bibinfo{pages}{253--264}.
\bibitem[{Liu et~al.(2017)Liu, Liu, and Peng}]{Liu2017}
\bibinfo{author}{Z.~Liu}, \bibinfo{author}{Z.~Liu}, \bibinfo{author}{Y.~Peng},
\newblock \bibinfo{title}{{Dimension reduction of Karhunen-Loeve expansion for
  simulation of stochastic processes}},
\newblock \bibinfo{journal}{Journal of Sound and Vibration}
  \bibinfo{volume}{408} (\bibinfo{year}{2017}) \bibinfo{pages}{168--189}.
\bibitem[{Grigoriu(1998)}]{Grigoriu1998}
\bibinfo{author}{M.~Grigoriu},
\newblock \bibinfo{title}{{Simulation of Stationary Non-Gaussian Translation
  Processes}},
\newblock \bibinfo{journal}{Journal of Engineering Mechanics}
  \bibinfo{volume}{124} (\bibinfo{year}{1998}) \bibinfo{pages}{121--126}.
\bibitem[{Shields et~al.(2011)Shields, Deodatis, and Bocchini}]{Shields2011}
\bibinfo{author}{M.~D. Shields}, \bibinfo{author}{G.~Deodatis},
  \bibinfo{author}{P.~Bocchini},
\newblock \bibinfo{title}{{A simple and efficient methodology to approximate a
  general non-Gaussian stationary stochastic process by a translation
  process}},
\newblock \bibinfo{journal}{Probabilistic Engineering Mechanics}
  \bibinfo{volume}{26} (\bibinfo{year}{2011}) \bibinfo{pages}{511--519}.
\bibitem[{Shields and Deodatis(2013)}]{Shields2013}
\bibinfo{author}{M.~D. Shields}, \bibinfo{author}{G.~Deodatis},
\newblock \bibinfo{title}{{A simple and efficient methodology to approximate a
  general non-Gaussian stationary stochastic vector process by a translation
  process with applications in wind velocity simulation}},
\newblock \bibinfo{journal}{Probabilistic Engineering Mechanics}
  \bibinfo{volume}{31} (\bibinfo{year}{2013}) \bibinfo{pages}{19--29}.
\bibitem[{Kim and Shields(2015)}]{Kim2015}
\bibinfo{author}{H.~Kim}, \bibinfo{author}{M.~D. Shields},
\newblock \bibinfo{title}{{Modeling strongly non-Gaussian non-stationary
  stochastic processes using the Iterative Translation Approximation Method and
  Karhunen-Lo{\`{e}}ve expansion}},
\newblock \bibinfo{journal}{Computers and Structures} \bibinfo{volume}{161}
  (\bibinfo{year}{2015}) \bibinfo{pages}{31--42}.
\bibitem[{Shields and Kim(2017)}]{Shields2017}
\bibinfo{author}{M.~D. Shields}, \bibinfo{author}{H.~Kim},
\newblock \bibinfo{title}{{Simulation of higher-order stochastic processes by
  spectral representation}},
\newblock \bibinfo{journal}{Probabilistic Engineering Mechanics}
  \bibinfo{volume}{47} (\bibinfo{year}{2017}) \bibinfo{pages}{1--15}.
\bibitem[{Vandanapu and Shields(2021)}]{Vandanapu2021}
\bibinfo{author}{L.~Vandanapu}, \bibinfo{author}{M.~D. Shields},
\newblock \bibinfo{title}{{3rd-order Spectral Representation Method: Simulation
  of multi-dimensional random fields and ergodic multi-variate random processes
  with fast Fourier transform implementation}},
\newblock \bibinfo{journal}{Probabilistic Engineering Mechanics}
  \bibinfo{volume}{64} (\bibinfo{year}{2021}) \bibinfo{pages}{103128}.
\bibitem[{Priestley(1965)}]{Priestley1965}
\bibinfo{author}{M.~B. Priestley},
\newblock \bibinfo{title}{Evolutionary spectra and non-stationary processes},
\newblock \bibinfo{journal}{Journal of the Royal Statistical Society. Series B
  (Methodological)} \bibinfo{volume}{27} (\bibinfo{year}{1965})
  \bibinfo{pages}{204--237}.
\bibitem[{Li and Kareem(1997)}]{Li1997}
\bibinfo{author}{Y.~Li}, \bibinfo{author}{A.~Kareem},
\newblock \bibinfo{title}{Simulation of {Multivariate} {Nonstationary} {Random}
  {Processes}: {Hybrid} {DFT} and {Digital} {Filtering} {Approach}},
\newblock \bibinfo{journal}{Journal of Engineering Mechanics}
  \bibinfo{volume}{123} (\bibinfo{year}{1997}) \bibinfo{pages}{1302--1310}.
\bibitem[{Huang(2014)}]{Huang2014}
\bibinfo{author}{G.~Huang},
\newblock \bibinfo{title}{{An efficient simulation approach for multivariate
  nonstationary process: Hybrid of wavelet and spectral representation
  method}},
\newblock \bibinfo{journal}{Probabilistic Engineering Mechanics}
  \bibinfo{volume}{37} (\bibinfo{year}{2014}) \bibinfo{pages}{74--83}.
\bibitem[{Li and Kareem(1991)}]{Li1991}
\bibinfo{author}{Y.~Li}, \bibinfo{author}{A.~Kareem},
\newblock \bibinfo{title}{{Simulation of Multivariate Nonstationary Random
  Processes by FFT}},
\newblock \bibinfo{journal}{Journal of Engineering Mechanics}
  \bibinfo{volume}{117} (\bibinfo{year}{1991}) \bibinfo{pages}{1037--1058}.
\bibitem[{Huang(2015)}]{Huang2015}
\bibinfo{author}{G.~Huang},
\newblock \bibinfo{title}{{Application of Proper Orthogonal Decomposition in
  Fast Fourier Transform—Assisted Multivariate Nonstationary Process
  Simulation}},
\newblock \bibinfo{journal}{Journal of Engineering Mechanics}
  \bibinfo{volume}{141} (\bibinfo{year}{2015}) \bibinfo{pages}{04015015}.
\bibitem[{Peng et~al.(2017)Peng, Huang, Chen, and Kareem}]{Peng2017}
\bibinfo{author}{L.~Peng}, \bibinfo{author}{G.~Huang},
  \bibinfo{author}{X.~Chen}, \bibinfo{author}{A.~Kareem},
\newblock \bibinfo{title}{{Simulation of multivariate nonstationary random
  processes: Hybrid stochastic wave and proper orthogonal decomposition
  approach}},
\newblock \bibinfo{journal}{Journal of Engineering Mechanics}
  \bibinfo{volume}{143} (\bibinfo{year}{2017}) \bibinfo{pages}{1--16}.
\bibitem[{Sakamoto and Ghanem(2002)}]{sakamoto2002simulation}
\bibinfo{author}{S.~Sakamoto}, \bibinfo{author}{R.~Ghanem},
\newblock \bibinfo{title}{Simulation of multi-dimensional non-gaussian
  non-stationary random fields},
\newblock \bibinfo{journal}{Probabilistic Engineering Mechanics}
  \bibinfo{volume}{17} (\bibinfo{year}{2002}) \bibinfo{pages}{167--176}.
\bibitem[{Ferrante and Graham-Brady(2005)}]{ferrante2005stochastic}
\bibinfo{author}{F.~Ferrante}, \bibinfo{author}{L.~Graham-Brady},
\newblock \bibinfo{title}{Stochastic simulation of non-gaussian/non-stationary
  properties in a functionally graded plate},
\newblock \bibinfo{journal}{Computer Methods in Applied Mechanics and
  Engineering} \bibinfo{volume}{194} (\bibinfo{year}{2005})
  \bibinfo{pages}{1675--1692}.
\bibitem[{Shields and Deodatis(2013)}]{shields2013estimation}
\bibinfo{author}{M.~Shields}, \bibinfo{author}{G.~Deodatis},
\newblock \bibinfo{title}{Estimation of evolutionary spectra for simulation of
  non-stationary and non-gaussian stochastic processes},
\newblock \bibinfo{journal}{Computers \& Structures} \bibinfo{volume}{126}
  (\bibinfo{year}{2013}) \bibinfo{pages}{149--163}.
\bibitem[{Dai et~al.(2019)Dai, Zheng, and Ma}]{dai2019explicit}
\bibinfo{author}{H.~Dai}, \bibinfo{author}{Z.~Zheng}, \bibinfo{author}{H.~Ma},
\newblock \bibinfo{title}{An explicit method for simulating non-gaussian and
  non-stationary stochastic processes by karhunen-lo{\`e}ve and polynomial
  chaos expansion},
\newblock \bibinfo{journal}{Mechanical Systems and Signal Processing}
  \bibinfo{volume}{115} (\bibinfo{year}{2019}) \bibinfo{pages}{1--13}.
\bibitem[{Montoya-Noguera et~al.(2019)Montoya-Noguera, Zhao, Hu, Wang, and
  Phoon}]{montoya2019simulation}
\bibinfo{author}{S.~Montoya-Noguera}, \bibinfo{author}{T.~Zhao},
  \bibinfo{author}{Y.~Hu}, \bibinfo{author}{Y.~Wang}, \bibinfo{author}{K.-K.
  Phoon},
\newblock \bibinfo{title}{Simulation of non-stationary non-gaussian random
  fields from sparse measurements using bayesian compressive sampling and
  karhunen-lo{\`e}ve expansion},
\newblock \bibinfo{journal}{Structural Safety} \bibinfo{volume}{79}
  (\bibinfo{year}{2019}) \bibinfo{pages}{66--79}.
\bibitem[{Zheng et~al.(2021)Zheng, Dai, Wang, and Wang}]{zheng2021sample}
\bibinfo{author}{Z.~Zheng}, \bibinfo{author}{H.~Dai},
  \bibinfo{author}{Y.~Wang}, \bibinfo{author}{W.~Wang},
\newblock \bibinfo{title}{A sample-based iterative scheme for simulating
  non-stationary non-gaussian stochastic processes},
\newblock \bibinfo{journal}{Mechanical Systems and Signal Processing}
  \bibinfo{volume}{151} (\bibinfo{year}{2021}) \bibinfo{pages}{107420}.
\bibitem[{Cramer(1967)}]{cramer1967}
\bibinfo{author}{H.~Cramer}, \bibinfo{title}{Stationary and related stochastic
  processes; sample function properties and their applications},
  \bibinfo{publisher}{Wiley}, \bibinfo{address}{New York},
  \bibinfo{year}{1967}.
\bibitem[{Priestley(1967)}]{Priestley1967}
\bibinfo{author}{M.~Priestley},
\newblock \bibinfo{title}{{Power spectral analysis of non-stationary random
  processes}},
\newblock \bibinfo{journal}{Journal of Sound and Vibration} \bibinfo{volume}{6}
  (\bibinfo{year}{1967}) \bibinfo{pages}{86--97}.
\bibitem[{Kolda and Bader(2009)}]{Kolda2009}
\bibinfo{author}{T.~G. Kolda}, \bibinfo{author}{B.~W. Bader},
\newblock \bibinfo{title}{{Tensor decompositions and applications}},
\newblock \bibinfo{journal}{SIAM Review} \bibinfo{volume}{51}
  (\bibinfo{year}{2009}) \bibinfo{pages}{455--500}.
\bibitem[{Clough and Penzien(1975)}]{clough_dynamics_1975}
\bibinfo{author}{R.~W. Clough}, \bibinfo{author}{J.~Penzien},
  \bibinfo{title}{Dynamics of structures}, \bibinfo{number}{xxii, 634 p.},
  \bibinfo{publisher}{McGraw-Hill}, \bibinfo{address}{New York},
  \bibinfo{year}{1975}.
\bibitem[{Liang et~al.(2007)Liang, Chaudhuri, and Shinozuka}]{Liang2007}
\bibinfo{author}{J.~Liang}, \bibinfo{author}{S.~R. Chaudhuri},
  \bibinfo{author}{M.~Shinozuka},
\newblock \bibinfo{title}{{Simulation of Nonstationary Stochastic Processes by
  Spectral Representation}},
\newblock \bibinfo{journal}{Journal of Engineering Mechanics}
  \bibinfo{volume}{133} (\bibinfo{year}{2007}) \bibinfo{pages}{616--627}.

\end{thebibliography}

\appendix

\section{Spectral Properties of the Modulated Orthogonal Increments}
\label{sec:AppendixA}





Herein, we demonstrate the orthogonality properties of the modulated orthogonal increments in the third-order non-stationary spectral representation given in Eqs. \eqref{eqn:ortho_U} and \eqref{eqn:ortho_V}.

Let us begin with the first-order orthogonality condition where it follows that
\begin{equation}
    \mathbb{E}[dU_{t}(\omega)] = \alpha(t, \omega)\mathbb{E}[dU(\omega)] + \beta(t, \omega)\mathbb{E}[dV(\omega)] = 0
\end{equation}
\begin{equation}
    \mathbb{E}[dV_{t}(\omega)] = \beta(t, \omega)\mathbb{E}[dU(\omega)] - \alpha(t, \omega)\mathbb{E}[dV(\omega)] = 0
\end{equation}
due to the fact that $dU(\omega)=dV(\omega)=0$

Next, consider the second-order orthogonality condition where
\begin{equation}
\begin{aligned}
    \mathbb{E}[dU_{t}^{2}(\omega)] &= \alpha^{2}(t, \omega)\mathbb{E}[dU^{2}(\omega)] + \beta^{2}(t, \omega)\mathbb{E}[dV^{2}(\omega)]\\
    & + 2\alpha(t, \omega)\beta(t, \omega)\mathbb{E}[dU(\omega)dV(\omega)]\\
    & = 2(\alpha^{2}(t, \omega) + \beta^{2}(t, \omega))S(\omega)d\omega\\
    & = 2|A(t, \omega)|^{2}S(\omega)d\omega\\
    & = 2S(t, \omega)d\omega
\end{aligned}
\end{equation}
\begin{equation}
\begin{aligned}
    \mathbb{E}[dV_{t}^{2}(\omega)] &= \beta^{2}(t, \omega)\mathbb{E}[dU^{2}(\omega)] + \alpha^{2}(t, \omega)\mathbb{E}[dV^{2}(\omega)]\\
    & - 2\alpha(t, \omega)\beta(t, \omega)\mathbb{E}[dU(\omega)dV(\omega)]\\
    & = 2(\alpha^{2}(t, \omega) + \beta^{2}(t, \omega))S(\omega)d\omega\\
    & = 2|A(t, \omega)|^{2}S(\omega)d\omega\\
    & = 2S(t, \omega)d\omega
\end{aligned}
\end{equation}

Finally, the third-order orthogonality condition yields
\begin{equation}
\begin{aligned}
    & \mathbb{E}[dU_{t}(\omega_{1})dU_{t}(\omega_{2})dU_{t}(\omega_{1} + \omega_{2})]\\
    & = \alpha(t, \omega_{1})\alpha(t, \omega_{2})\alpha(t, \omega_{1} + \omega_{2})\mathbb{E}[dU(\omega_{1})dU(\omega_{2})dU(\omega_{1} + \omega_{2})]\\
    & + \alpha(t, \omega_{1})\alpha(t, \omega_{2})\beta(t, \omega_{1} + \omega_{2})\mathbb{E}[dU(\omega_{1})dU(\omega_{2})dV(\omega_{1} + \omega_{2})]\\
    & + \alpha(t, \omega_{1})\beta(t, \omega_{2})\alpha(t, \omega_{1} + \omega_{2})\mathbb{E}[dU(\omega_{1})dV(\omega_{2})dU(\omega_{1} + \omega_{2})]\\
    & + \beta(t, \omega_{1})\alpha(t, \omega_{2})\alpha(t, \omega_{1} + \omega_{2})\mathbb{E}[dV(\omega_{1})dU(\omega_{2})dU(\omega_{1} + \omega_{2})]\\
    & + \alpha(t, \omega_{1})\beta(t, \omega_{2})\beta(t, \omega_{1} + \omega_{2})\mathbb{E}[dU(\omega_{1})dV(\omega_{2})dV(\omega_{1} + \omega_{2})]\\
    & + \beta(t, \omega_{1})\alpha(t, \omega_{2})\beta(t, \omega_{1} + \omega_{2})\mathbb{E}[dV(\omega_{1})dU(\omega_{2})dV(\omega_{1} + \omega_{2})]\\
    & + \beta(t, \omega_{1})\beta(t, \omega_{2})\alpha(t, \omega_{1} + \omega_{2})\mathbb{E}[dV(\omega_{1})dV(\omega_{2})dU(\omega_{1} + \omega_{2})]\\
    & + \beta(t, \omega_{1})\beta(t, \omega_{2})\beta(t, \omega_{1} + \omega_{2})\mathbb{E}[dV(\omega_{1})dV(\omega_{2})dV(\omega_{1} + \omega_{2})]\\
    & = 2\alpha(t, \omega_{1})\alpha(t, \omega_{2})\alpha(t, \omega_{1} + \omega_{2})\mathbb{R}B(\omega_{1}, \omega_{2})d\omega_{1}d\omega_{2}\\
    & - 2 \alpha(t, \omega_{1})\alpha(t, \omega_{2})\beta(t, \omega_{1} + \omega_{2})\mathbb{I}B(\omega_{1}, \omega_{2})d\omega_{1}d\omega_{2}\\
    & + 2 \alpha(t, \omega_{1})\beta(t, \omega_{2})\alpha(t, \omega_{1} + \omega_{2})\mathbb{I}B(\omega_{1}, \omega_{2})d\omega_{1}d\omega_{2}\\
    & + 2 \beta(t, \omega_{1})\alpha(t, \omega_{2})\alpha(t, \omega_{1} + \omega_{2})\mathbb{I}B(\omega_{1}, \omega_{2})d\omega_{1}d\omega_{2}\\
    & + 2 \alpha(t, \omega_{1})\beta(t, \omega_{2})\beta(t, \omega_{1} + \omega_{2})\mathbb{R}B(\omega_{1}, \omega_{2})d\omega_{1}d\omega_{2}\\
    & + 2 \beta(t, \omega_{1})\alpha(t, \omega_{2})\beta(t, \omega_{1} + \omega_{2})\mathbb{R}B(\omega_{1}, \omega_{2})d\omega_{1}d\omega_{2}\\
    & - 2 \beta(t, \omega_{1})\beta(t, \omega_{2})\alpha(t, \omega_{1} + \omega_{2})\mathbb{R}B(\omega_{1}, \omega_{2})d\omega_{1}d\omega_{2}\\
    & + 2 \beta(t, \omega_{1})\beta(t, \omega_{2})\beta(t, \omega_{1} + \omega_{2})\mathbb{I}B(\omega_{1}, \omega_{2})d\omega_{1}d\omega_{2}\\
    & = 2A^{*}(t, \omega_{1})A^{*}(t, \omega_{2})A(t, \omega_{1} + \omega_{2})\mathbb{E}[dZ^{*}(\omega_{1})dZ^{*}(\omega_{2})dZ^{*}(\omega_{1} + \omega_{2})]\\
    & = 2B(t, \omega_{1}, \omega_{2})d\omega_{1}d\omega_{2}\\
\end{aligned}
\end{equation}

\begin{equation}
\begin{aligned}
    & \mathbb{E}[dV_{t}(\omega_{1})dV_{t}(\omega_{2})dV_{t}(\omega_{1} + \omega_{2})]\\
    & = \beta(t, \omega_{1})\beta(t, \omega_{2})\beta(t, \omega_{1} + \omega_{2})\mathbb{E}[dU(\omega_{1})dU(\omega_{2})dU(\omega_{1} + \omega_{2})]\\
    & - \beta(t, \omega_{1})\beta(t, \omega_{2})\alpha(t, \omega_{1} + \omega_{2})\mathbb{E}[dU(\omega_{1})dU(\omega_{2})dV(\omega_{1} + \omega_{2})]\\
    & - \beta(t, \omega_{1})\alpha(t, \omega_{2})\beta(t, \omega_{1} + \omega_{2})\mathbb{E}[dU(\omega_{1})dV(\omega_{2})dU(\omega_{1} + \omega_{2})]\\
    & - \alpha(t, \omega_{1})\beta(t, \omega_{2})\beta(t, \omega_{1} + \omega_{2})\mathbb{E}[dV(\omega_{1})dU(\omega_{2})dU(\omega_{1} + \omega_{2})]\\
    & + \beta(t, \omega_{1})\alpha(t, \omega_{2})\alpha(t, \omega_{1} + \omega_{2})\mathbb{E}[dU(\omega_{1})dV(\omega_{2})dV(\omega_{1} + \omega_{2})]\\
    & + \alpha(t, \omega_{1})\beta(t, \omega_{2})\alpha(t, \omega_{1} + \omega_{2})\mathbb{E}[dV(\omega_{1})dU(\omega_{2})dV(\omega_{1} + \omega_{2})]\\
    & + \alpha(t, \omega_{1})\alpha(t, \omega_{2})\beta(t, \omega_{1} + \omega_{2})\mathbb{E}[dV(\omega_{1})dV(\omega_{2})dU(\omega_{1} + \omega_{2})]\\
    & - \alpha(t, \omega_{1})\alpha(t, \omega_{2})\alpha(t, \omega_{1} + \omega_{2})\mathbb{E}[dV(\omega_{1})dV(\omega_{2})dV(\omega_{1} + \omega_{2})]\\
    & = 2\beta(t, \omega_{1})\beta(t, \omega_{2})\beta(t, \omega_{1} + \omega_{2})\mathbb{R}B(\omega_{1}, \omega_{2})d\omega_{1}d\omega_{2}\\
    & + 2 \beta(t, \omega_{1})\beta(t, \omega_{2})\alpha(t, \omega_{1} + \omega_{2})\mathbb{I}B(\omega_{1}, \omega_{2})d\omega_{1}d\omega_{2}\\
    & - 2 \beta(t, \omega_{1})\alpha(t, \omega_{2})\beta(t, \omega_{1} + \omega_{2})\mathbb{I}B(\omega_{1}, \omega_{2})d\omega_{1}d\omega_{2}\\
    & - 2 \alpha(t, \omega_{1})\beta(t, \omega_{2})\beta(t, \omega_{1} + \omega_{2})\mathbb{I}B(\omega_{1}, \omega_{2})d\omega_{1}d\omega_{2}\\
    & + 2 \beta(t, \omega_{1})\alpha(t, \omega_{2})\alpha(t, \omega_{1} + \omega_{2})\mathbb{R}B(\omega_{1}, \omega_{2})d\omega_{1}d\omega_{2}\\
    & + 2 \alpha(t, \omega_{1})\beta(t, \omega_{2})\alpha(t, \omega_{1} + \omega_{2})\mathbb{R}B(\omega_{1}, \omega_{2})d\omega_{1}d\omega_{2}\\
    & - 2 \alpha(t, \omega_{1})\alpha(t, \omega_{2})\beta(t, \omega_{1} + \omega_{2})\mathbb{R}B(\omega_{1}, \omega_{2})d\omega_{1}d\omega_{2}\\
    & - 2 \alpha(t, \omega_{1})\alpha(t, \omega_{2})\alpha(t, \omega_{1} + \omega_{2})\mathbb{I}B(\omega_{1}, \omega_{2})d\omega_{1}d\omega_{2}\\
    & = 2A^{*}(t, \omega_{1})A^{*}(t, \omega_{2})A(t, \omega_{1} + \omega_{2})\mathbb{E}[dZ^{*}(\omega_{1})dZ^{*}(\omega_{2})dZ^{*}(\omega_{1} + \omega_{2})]\\
    & = 2B(t, \omega_{1}, \omega_{2})d\omega_{1}d\omega_{2}\\
\end{aligned}
\end{equation}


\section{Ensemble Properties}
\label{ap:B}

Here, we show that the expansion in Eq.\ \eqref{eqn:process_infinite_sum} satisfies the necessary ensemble properties in 1st-order (mean), 2nd-order (variance) and 3rd-order respectively

\subsection{Mean}

The processes are assumed to be of zero-mean, $\mathbb{E}[X(t)] = 0$. Since expectation is commutative over summation
\begin{equation}
\begin{aligned}
    \mathbb{E}[X(t)] &=  \mathbb{E}[\sqrt{2}\sum_{k=0}^{N-1}[2S_{p}(t, \omega_{k})\Delta\omega]^{\frac{1}{2}}\cos(\omega_{k}t + \phi_{k})\\
    &+ \sqrt{2}\sum_{k=0}^{N-1}\sum_{i+j=k}^{i \geq j \geq 0}[2S(t, \omega_{k})\Delta\omega]^{\frac{1}{2}}b_{p}(t, \omega_{i}, \omega_{j}) cos(\omega_{k}t + \phi_{i} + \phi_{j} + \beta(t, \phi_{i}, \phi_{j})]\\
    &=  \sqrt{2}\sum_{k=0}^{N-1}[2S_{p}(t, \omega_{k})\Delta\omega]^{\frac{1}{2}}\mathbb{E}[cos(\omega_{k}t + \phi_{k})]\\
    &+ \sqrt{2}\sum_{k=0}^{N-1}\sum_{i+j=k}^{i \geq j \geq 0}[2S(t, \omega_{k})\Delta\omega]^{\frac{1}{2}}b_{p}(t, \omega_{i}, \omega_{j})\mathbb{E}[cos(\omega_{k}t + \phi_{i} + \phi_{j} + \beta(t, \phi_{i}, \phi_{j})]\\
\end{aligned}
\end{equation}
Computing the expectation of $\cos(\omega_{k}t + \phi_{k})$ and $\cos(\omega_{k}t + \phi_{i} + \phi_{j} + \beta(t, \omega_{i}, \omega_{j}))$ we have
\begin{equation}
    \mathbb{E}[\cos(\omega_{k}t + \phi_{k})] = \int_{-\infty}^{\infty}p_{\phi_{k}}\cos(\omega_{k}t + \phi_{k})d\phi_{k} = \frac{1}{2\pi}\int_{0}^{2\pi}\cos(\omega_{k}t + \phi_{k})d\phi_{k} = 0
\end{equation}
\begin{equation}
\begin{aligned}
    \mathbb{E}[\cos(\omega_{k}t + \phi_{i} + \phi_{j} + \beta(t, \omega_{i}, \omega_{j}))] 
    &= \int_{-\infty}^{\infty}\int_{-\infty}^{\infty}p_{\phi_{i}}p_{\phi_{j}}\cos(\omega_{k}t + \phi_{i} + \phi_{j} + \beta(t, \omega_{i}, \omega_{j}))d\phi_{i}d\phi_{j} \\
    &= \frac{1}{4\pi^{2}}\int_{0}^{2\pi}\int_{0}^{2\pi}\cos(\omega_{k}t + \phi_{i} + \phi_{j} + \beta(t, \omega_{i}, \omega_{j}))d\phi_{i}d\phi_{j} = 0
\end{aligned}
\end{equation}
Using the above results, we can conclude that
\begin{equation}
    \mathbb{E}[X(t)] = 0
\end{equation}

\subsection{2-point Correlation Function}

The $2$-point correlation function can be computed as
\begin{equation}
\begin{aligned}
    &R_{2}(t, t + \tau) = \mathbb{E}[X(t)X(t + \tau)]\\
    &= \mathbb{E}[(\sqrt{2}\sum_{c=0}^{N-1}[2S_{p}(t, \omega_{c})\Delta\omega]^{\frac{1}{2}}\cos(\omega_{c}t + \phi_{c})\\
    &+ \sqrt{2}\sum_{c=0}^{N-1}\sum_{a+b=c}^{a \geq b \geq 0}[2S(t, \omega_{c})\Delta\omega]^{\frac{1}{2}}b_{p}(t, \omega_{a}, \omega_{b})\cos(\omega_{c}t + \phi_{a} + \phi_{b} + \beta(t, \phi_{a}, \phi_{b}))\\
    & (\sqrt{2}\sum_{k=0}^{N-1}[2S_{p}(t + \tau, \omega_{k})\Delta\omega]^{\frac{1}{2}}\cos(\omega_{k}(t + \tau) + \phi_{k})\\
    &+ \sqrt{2}\sum_{k=0}^{N-1}\sum_{i+j=k}^{i \geq j \geq 0}[2S(t + \tau, \omega_{k})\Delta\omega]^{\frac{1}{2}}b_{p}(t + \tau, \omega_{i}, \omega_{j})\cos(\omega_{k}(t + \tau) + \phi_{i} + \phi_{j} + \beta(t + \tau, \phi_{i}, \phi_{j})))]\\
    &= 4 \mathbb{E}[\sum_{c=0}^{N-1}\sum_{k=0}^{N-1} \sqrt{S_{p}(t, \omega_{c}) S_{p}(t + \tau, \omega_{k})} \Delta \omega \cos(\omega_{c}t + \phi_{c})\cos(\omega_{k}(t + \tau) + \phi_{k})\\
    &+ \sum_{k=0}^{N-1} \sum_{c=0}^{N-1}\sum_{a+b=c}^{a \geq b \geq 0} \sqrt{S(t + \tau, \omega_{k}) S(t, \omega_{c})} b_{p}(t, \omega_{a}, \omega_{b}) \Delta \omega \\
    & \cos(\omega_{k}(t + \tau) + \phi_{k})\cos(\omega_{c}t + \phi_{a} + \phi_{b} + \beta(t, \phi_{a}, \phi_{b}))\\
    &+ \sum_{c=0}^{N-1} \sum_{k=0}^{N-1}\sum_{i+j=k}^{i \geq j \geq 0} \sqrt{S_{p}(t, \omega_{c})S(t + \tau, \omega_{k})} b_{p}(t + \tau, \omega_{i}, \omega_{j}) \Delta \omega \\
    & \cos(\omega_{c}t + \phi_{c}) \cos(\omega_{k}(t + \tau) + \phi_{i} + \phi_{j} + \beta(t + \tau, \phi_{i}, \phi_{j})))\\
    & + \sum_{c=0}^{N-1} \sum_{k=0}^{N-1}\sum_{a+b=c}^{a \geq b \geq 0}\sum_{i+j=k}^{i \geq j \geq 0} \sqrt{S(t, \omega_{c})S(t + \tau, \omega_{k})} b_{p}(t, \omega_{a}, \omega_{b}) b_{p}(t + \tau, \omega_{i}, \omega_{j}) \Delta \omega \\
    & \cos(\omega_{c}t + \phi_{a} + \phi_{b} + \beta(t, \phi_{a}, \phi_{b})) \cos(\omega_{k}(t + \tau) + \phi_{i} + \phi_{j} + \beta(t + \tau, \phi_{i}, \phi_{j})) ]\\
    &= 2 [\sum_{k=0}^{N-1} \sqrt{S_{p}(t, \omega_{k}) S_{p}(t + \tau, \omega_{k})} \Delta \omega \cos(\omega_{k}\tau) \\
    & + \sum_{k=0}^{N-1}\sum_{i+j=k}^{i \geq j \geq 0} \sqrt{S(t, \omega_{k})S(t + \tau, \omega_{k})} b_{p}(t, \omega_{a}, \omega_{b}) b_{p}(t + \tau, \omega_{i}, \omega_{j}) \Delta \omega \cos(\omega_{k}\tau)] \\
    &= 2 \sum_{k=0}^{N-1} [ \sqrt{S_{p}(t, \omega_{k}) S_{p}(t + \tau, \omega_{k})} \Delta \omega \cos(\omega_{k}\tau) \\
    & + \sum_{i+j=k}^{i \geq j \geq 0} \sqrt{S(t, \omega_{k})S(t + \tau, \omega_{k})} b_{p}(t, \omega_{i}, \omega_{j})b_{p}(t + \tau, \omega_{i}, \omega_{j}) \Delta \omega \cos(\omega_{k}\tau + \beta(t + \tau, \phi_{i}, \phi_{j}) - \beta(t, \phi_{i}, \phi_{j}))] \\
\end{aligned}
\end{equation}
Substituting the value of $\tau=0$ i.e. computing the variance at time $t$,
\begin{equation}
\begin{aligned}
    & R(t, t) = E[X^2(t)] = 2 \sum_{k=0}^{N-1} S(t, \omega_{k}) \Delta\omega \\
    & R(t, t) = E[X^2(t)] = \int_{-\infty}^{\infty} S(t, \omega) d \omega\\
\end{aligned}
\end{equation}
we see that the process has variance equal to the integral of the evolutionary spectrum, as expected. This is equivalent to the variance function from the $2^{nd}$-order SRM \cite{Liang2007}.

\subsection{3-point Correlation Function}

The $3$-point correlation function can be computed as
\begin{equation}
\begin{aligned}
    &R_{3}(t, t + \tau_{1}, t + \tau_{2}) = \mathbb{E}[X(t)X(t + \tau_{1})X(t + \tau_{2})]\\
    &= \mathbb{E}[(\sqrt{2}\sum_{c=0}^{N-1}[2S_{p}(t, \omega_{c})\Delta\omega]^{\frac{1}{2}}\cos(\omega_{c}t + \phi_{c})\\
    &+ \sqrt{2}\sum_{c=0}^{N-1}\sum_{a+b=c}^{a \geq b \geq 0}[2S(t, \omega_{c})\Delta\omega]^{\frac{1}{2}}b_{p}(t, \omega_{a}, \omega_{b})\cos(\omega_{c}t + \phi_{a} + \phi_{b} + \beta(t, \phi_{a}, \phi_{b}))\\
    & (\sqrt{2}\sum_{k=0}^{N-1}[2S_{p}(t + \tau_1, \omega_{k})\Delta\omega]^{\frac{1}{2}}\cos(\omega_{k}(t + \tau_1) + \phi_{k})\\
    &+ \sqrt{2}\sum_{k=0}^{N-1}\sum_{i+j=k}^{i \geq j \geq 0}[2S(t + \tau_1, \omega_{k})\Delta\omega]^{\frac{1}{2}}b_{p}(t + \tau_1, \omega_{i}, \omega_{j})\cos(\omega_{k}(t + \tau_1) + \phi_{i} + \phi_{j} + \beta(t + \tau_1, \phi_{i}, \phi_{j}))) \\
    & (\sqrt{2}\sum_{z=0}^{N-1}[2S_{p}(t + \tau_2, \omega_{z})\Delta\omega]^{\frac{1}{2}}\cos(\omega_{z}(t + \tau_2) + \phi_{z})\\
    &+ \sqrt{2}\sum_{z=0}^{N-1}\sum_{x+y=z}^{x \geq y \geq 0}[2S(t + \tau_2, \omega_{z})\Delta\omega]^{\frac{1}{2}}b_{p}(t + \tau_2, \omega_{x}, \omega_{y})\cos(\omega_{z}(t + \tau_2) + \phi_{x} + \phi_{y} + \beta(t + \tau_2, \phi_{x}, \phi_{y}))) ]\\
\end{aligned}
\end{equation}
which simplifies to
\begin{equation}
\begin{aligned}
    &R_{3}(t, t + \tau_{1}, t + \tau_{2}) = \\
    & 2 (\Delta \omega)^{3/2} [ \sum_{k=0}^{N-1}\sum_{i+j=k}^{i \geq j \geq 0} \sqrt{S_{p}(t, \omega_{i})S_{p}(t + \tau_{1}, \omega_{j})S(t+ \tau_{2}, \omega_{k})}b_{p}(t + \tau_{2}, \omega_{i}, \omega_{j})\cos(\omega_{k}\tau_2 + \beta(t + \tau_2, \phi_{i}, \phi_{j}) - \omega_{j}\tau_1)\\
    & + \sqrt{S_{p}(t + \tau_{1}, \omega_{i})S_{p}(t, \omega_{j})S(t+ \tau_{2}, \omega_{k})}b_{p}(t + \tau_{2}, \omega_{i}, \omega_{j})\cos(\omega_{k}\tau_2  + \beta(t + \tau_2, \phi_{i}, \phi_{j}) - \omega_{i}\tau_1)\\
    & + \sqrt{S_{p}(t + \tau_{2}, \omega_{i})S_{p}(t, \omega_{j})S(t+ \tau_{1}, \omega_{k})}b_{p}(t + \tau_{1}, \omega_{i}, \omega_{j})\cos(\omega_{k}\tau_{1} + \beta(t + \tau_1, \phi_{i}, \phi_{j}) - \omega_{i}\tau_{2})\\
    & + \sqrt{S_{p}(t, \omega_{i})S_{p}(t + \tau_{2}, \omega_{j})S(t+ \tau_{1}, \omega_{k})}b_{p}(t + \tau_{1}, \omega_{i}, \omega_{j})\cos(\omega_{k}\tau_{1} + \beta(t + \tau_1, \phi_{i}, \phi_{j}) - \omega_{j}\tau_{2})\\
    & + \sqrt{S_{p}(t + \tau_{1}, \omega_{i})S_{p}(t + \tau_{2}, \omega_{j})S(t, \omega_{k})}b_{p}(t, \omega_{i}, \omega_{j})\cos(\beta(t, \phi_{i}, \phi_{j}) - \omega_{i}\tau_{1} - \omega_{j}\tau_{2})\\
    & + \sqrt{S_{p}(t + \tau_{2}, \omega_{i})S_{p}(t + \tau_{1}, \omega_{j})S(t, \omega_{k})}b_{p}(t, \omega_{i}, \omega_{j})\cos(\beta(t + \tau_2, \phi_{i}, \phi_{j}) - \omega_{i}\tau_{2} - \omega_{j}\tau_{1}) ] \\
\end{aligned}
\end{equation}
Substituting $\tau_1 = \tau_2 = 0$, i.e. computing the 3-order correlation function at time $t$,
\begin{equation}
\begin{aligned}
    & R_{3}(t, t, t) = 12 (\Delta \omega)^{3/2} \sum_{k=0}^{N-1}\sum_{i+j=k}^{i \geq j \geq 0} \sqrt{S_{p}(t, \omega_{i})S_{p}(t, \omega_{j})S(t, \omega_{k})}b_{p}(t, \omega_{i}, \omega_{j})\cos(\beta(t, \phi_{i}, \phi_{j})) \\
\end{aligned}
\end{equation}
Plugging in the value of $b_{p}(t, \omega_1, \omega_2)$, yields
\begin{equation}
\begin{aligned}
    & R_{3}(t, t, t) = E[X^3(t)] = 12 (\Delta \omega)^2 \sum_{k=0}^{N-1}\sum_{i+j=k}^{i \geq j \geq 0} B(t, \omega_{i}, \omega_{j}) \\
    & R_{3}(t, t, t) = E[X^3(t)] = \int_{-\infty}^{\infty} B(t, \omega_{i}, \omega_{j}) d\omega_1 d\omega_2,\\
\end{aligned}
\end{equation}
which is consistent with the fact that the integral of the evolutionary bispectrum is equal to the third moment of the process at time $t$ as shown in \ref{ap:C}.

\section{Evolutionary bispectrum as the distribution of skewness}
\label{ap:C}

Here, we show that the proposed evolutionary bispectrum, when integrated over frequency pairs $\omega_1, \omega_2$, yields the third-moment of the process at a given time $t$.

First, consider the three-point correlation given by
\begin{equation}
\begin{aligned}
    R_3(s, t, u) & = \mathbb{E}\left[\int_{-\infty}^{\infty}\int_{-\infty}^{\infty}\int_{-\infty}^{\infty} A(s, \omega_1)dZ(\omega_1)A(t, \omega_2)dZ(\omega_2)A(u, \omega_3)dZ(\omega_3)\right] \\
    & = \int_{-\infty}^{\infty}\int_{-\infty}^{\infty}\int_{-\infty}^{\infty} A(s, \omega_1)A(t, \omega_2)A(u, \omega_3)\mathbb{E}\left[dZ(\omega_1)dZ(\omega_2)dZ(\omega_3)\right] \\
\end{aligned}
\end{equation}
Next, consider the third-order orthogonality conditions of the spectral process, here repeated from Eq.\ \eqref{eqn:3_order_orthogonal} 
\begin{equation}
\begin{aligned}
    & \mathbb{E}[dZ(\omega_1)dZ(\omega_2)dZ(\omega_3)] = B(\omega_1, \omega_2)d\omega_1 d\omega_2 \quad &\text{if} \quad  \omega_1 + \omega_2 - \omega_3 = 0\\
    & \mathbb{E}[dZ(\omega_1)dZ(\omega_2)dZ(\omega_3)] = 0 & \text{otherwise} \\
\end{aligned}
\end{equation}
Applying these orthogonality conditions and considering the third moment of a point at time $t$ yields
\begin{equation}
\begin{aligned}
    & R_3(s, t, u) = \int_{-\infty}^{\infty}\int_{-\infty}^{\infty} A(s, \omega_1)A(t, \omega_2)A(u, \omega_1 + \omega_2)\mathbb{E}[dZ(\omega_1)dZ(\omega_2)dZ^{*}(\omega_1 + \omega_2)] \\
    & R_3(t, t, t) = E[X^3(t)]= \int_{-\infty}^{\infty}\int_{-\infty}^{\infty} A(t, \omega_1)A(t, \omega_2)A(t, \omega_1 + \omega_2)B(\omega_1, \omega_2)d\omega_1 d\omega_2 \\
\end{aligned}
\end{equation}
Finally, we recognize that the integrand is the proposed evolutionary bispectrum such that
\begin{equation}
    E[X^3(t)]= \int_{-\infty}^{\infty}\int_{-\infty}^{\infty} B(t,\omega_1,\omega_2)d\omega_1d\omega_2,
\end{equation}
which is related to the skewness, $\gamma_3$ through
\begin{equation}
    \gamma_3 = \frac{\mathbb{E}[X^{3}(t)]}{(\mathbb{E}[X^2(t)])^{3/2}}
\end{equation}


\section{Theoretical reconstruction of decomposed spectral quantities}
\label{ap:D}

\subsection{Pure component of the Variance}

From the orthogonal projection in Eq.\ \eqref{eqn:second_order_projection}], we have
\begin{equation}
    \sqrt{S_{p}(t, \omega)} = \sum_{q = 1}^{N_q} a_{q}(t)\Phi_q(\omega)
\end{equation}
The square of the above equation yields
\begin{equation}
    S_{p}(t, \omega) = \sum_{q_1 = 1}^{N_q} \sum_{q_2 = 1}^{N_q} a_{q_{1}}(t)\Phi_{q_{1}}(\omega) a_{q_{2}}(t)\Phi_{q_{2}}(\omega)
\end{equation}
Taking the sum of along the range of the frequency terms yields
\begin{equation}
    \sum_{k = 0}^{N_{\omega}-1} S_{p}(t, \omega_k) = \sum_{k = 0}^{N_{\omega} - 1} \sum_{q_1 = 1}^{N_q} \sum_{q_2 = 1}^{N_q} a_{q_{1}}(t)\Phi_{q_{1}}(\omega_k) a_{q_{2}}(t)\Phi_{q_{2}}(\omega_k)
    \label{eqn:apd1}
\end{equation}
Since $\Phi_{q}(\omega)$ are orthogonal vectors, the following equation holds
\begin{equation}
    \sum_{k = 1}^{N_{\omega}} \Phi_{q_{1}}(\omega_k) \Phi_{q_{2}}(\omega_k) = \begin{cases}
    1, \text{if } q_1 = q_2\\
    0, \text{otherwise}
\end{cases}
\end{equation}
    Substituting the above orthogonality condition into Eq.\ \eqref{eqn:apd1} yields
\begin{equation}
\begin{aligned}
    & \sum_{k = 1}^{N_{\omega}} S_{p}(t, \omega_k) = \sum_{q = 1}^{N_q} a_{q}^{2}(t) \sum_{k = 1}^{N_{\omega}} \Phi_{q}^{2}(\omega_k) \\
    & \sum_{k = 1}^{N_{\omega}} S_{p}(t, \omega_k) = \sum_{q = 1}^{N_q} a_{q}^{2}(t) \\
\end{aligned}
\end{equation}
The variance of the pure component of the non-stationary process $X(t)$ can then be written as
\begin{equation}
    Var[X_{p}(t)] = 2 \Delta \omega \sum_{k = 1}^{N_{\omega}} S_{p}(t, \omega_k) = 2 \Delta \omega \sum_{q = 1}^{N_q} a_{q}^{2}(t),
    \label{eqn:pure_variance_pod_theoretical}
\end{equation}
which illustrates that $a_q(t)$ can be interpreted as modulating the pure component of the variance.

\subsection{Interactive component of the Variance}

Consider the decomposition of the bispectrum from Eq.\ \eqref{eqn:bispectrum_approximation}
\begin{equation}
    \frac{|B(t, \omega_{i}, \omega_{j})|}{\sqrt{S_{p}(t, \omega_i)S_p(t, \omega_j)}} = \sum_{r=1, s=1}^{N_q, N_q} |b_{rs}(t)| \theta_{rs}(\omega_i, \omega_j) = \sum_{r=1, s=1}^{N_q, N_q} b_{rs}(t)\Phi_r(\omega_i)\Phi_{s}(\omega_j) 
\end{equation}
The square of the above equation yields
\begin{equation}
    \frac{|B(t, \omega_{i}, \omega_{j})|^2}{S_{p}(t, \omega_i)S_p(t, \omega_j)} = \sum_{r_1=1, s_1=1}^{N_q, N_q}\sum_{r_2=1, s_2=1}^{N_q, N_q} |b_{r_1s_1}(t)| \theta_{r_1s_1}(\omega_i, \omega_j) |b_{r_2s_2}(t)| \theta_{r_2s_2}(\omega_i, \omega_j)
\end{equation}
Taking the sums of the term along the frequency range yields
\begin{equation}
\begin{aligned}
    & \sum_{k = 1}^{N_\omega} \sum_{i+j=k}^{i \geq 0, j \geq 0}\frac{|B(t, \omega_{i}, \omega_{j})|^2}{S_{p}(t, \omega_i)S_p(t, \omega_j)} = \sum_{i+j=k}^{i \geq 0, j \geq 0} \sum_{r_1=1, s_1=1}^{N_q, N_q}\sum_{r_2=1, s_2=1}^{N_q, N_q} |b_{r_1s_1}(t)| \theta_{r_1s_1}(\omega_i, \omega_j) |b_{r_2s_2}(t)| \theta_{r_2s_2}(\omega_i, \omega_j) \\
    & = \sum_{k = 1}^{N_\omega} \sum_{i+j=k}^{i \geq 0, j \geq 0} \sum_{r_1=1, s_1=1}^{N_q, N_q}\sum_{r_2=1, s_2=1}^{N_q, N_q} |b_{r_1s_1}(t)| \Phi_{r_1}(\omega_i)\Phi_{s_1}(\omega_j) |b_{r_2s_2}(t)| \Phi_{r_2}(\omega_i)\Phi_{s_2}(\omega_j) \\
    & = \sum_{r_1=1, s_1=1}^{N_q, N_q}\sum_{r_2=1, s_2=1}^{N_q, N_q} |b_{r_1s_1}(t)||b_{r_2s_2}(t)| \sum_{k = 1}^{N_\omega} \sum_{i+j=k}^{i \geq 0, j \geq 0} \Phi_{r_1}(\omega_i)\Phi_{s_1}(\omega_j)\Phi_{r_2}(\omega_i)\Phi_{s_2}(\omega_j) \\
    & = \sum_{r_1=1, s_1=1}^{N_q, N_q}\sum_{r_2=1, s_2=1}^{N_q, N_q} |b_{r_1s_1}(t)| |b_{r_2s_2}(t)| \sum_{k = 1}^{N_\omega} \sum_{i+j=k}^{i \geq 0, j \geq 0} \Phi_{r_1}(\omega_i)\Phi_{r_2}(\omega_i)\Phi_{s_1}(\omega_j)\Phi_{s_2}(\omega_j) \\
    \label{eqn:D9}
\end{aligned}
\end{equation}
We have the orthogonality condition of the singular vectors
\begin{equation}
    \sum_{i=1, j=1}^{N_\omega, N_\omega} \Phi_{r_1}(\omega_i)\Phi_{r_2}(\omega_i)\Phi_{s_1}(\omega_j)\Phi_{s_2}(\omega_j) = \begin{cases}
    1, \text{if } r_1 = r_2, s_1 = s_2\\
    0, \text{otherwise}
    \end{cases}
\end{equation}
Assuming that most of the information in the singular vectors is contained in the lower frequency ranges, the following approximation holds
\begin{equation}
    \sum_{k = 1}^{N_\omega} \sum_{i + j = k}^{i \geq 0, j \geq 0} \Phi_{r_1}(\omega_i)\Phi_{r_2}(\omega_i)\Phi_{s_1}(\omega_j)\Phi_{s_2}(\omega_j) \simeq \begin{cases}
    1, \text{if } r_1 = r_2, s_1 = s_2\\
    0, \text{otherwise}
    \end{cases}
\end{equation}
Substituting this orthogonality condition into Eq.\ \eqref{eqn:D9} yields
\begin{equation}
    \sum_{k = 1}^{N_\omega} \sum_{i+j=k}^{i \geq 0, j \geq 0}\frac{|B(t, \omega_{i}, \omega_{j})|^2}{S_{p}(t, \omega_i)S_p(t, \omega_j)} = \sum_{r = 1, s = 1}^{N_q, N_q} |b_{rs}(t)|^{2}
\end{equation}
The variance of the interactive component of the non-stationary process $X(t)$ can be written as
\begin{equation}
    Var[X_{I}(t)] = \Delta \omega^{2} \sum_{k = 1}^{N_\omega} \sum_{i+j=k}^{i \geq 0, j \geq 0}\frac{|B(t, \omega_{i}, \omega_{j})|^2}{S_{p}(t, \omega_i)S_p(t, \omega_j)} = \Delta\omega^2 \sum_{r = 1, s = 1}^{N_q, N_q} |b_{rs}(t)|^{2},
    \label{eqn:interactive_variance_pod_theoretical}
\end{equation}
which confirms that $b_{rs}(t)$ modulates the variance associated with wave interactions.

\subsection{Third Moment}

The bispectrum can be represented as
\begin{equation}
\begin{aligned}
    & B(t, \omega_i, \omega_j) = \frac{B(t, \omega_{i}, \omega_{j})}{\sqrt{S_{p}(t, \omega_i)S_p(t, \omega_j)}} \sqrt{S_{p}(t, \omega_i)}\sqrt{S_p(t, \omega_j)} \\
    & = \sum_{r=1, s=1}^{N_q, N_q} b_{rs}(t)\theta_{rs}(\omega_i, \omega_j) \sum_{q_1 = 1}^{N_{q_1}} a_{q_{1}}(t)\Phi_{q_1}(\omega_i) \sum_{q_2 = 1}^{N_{q_2}} a_{q_{2}}(t)\Phi_{q_2}(\omega_j)\\
    & = \sum_{r=1, s=1}^{N_q, N_q} \sum_{q_1 = 1}^{N_q} \sum_{q_2 = 1}^{N_q} b_{rs}(t) a_{q_{1}}(t) a_{q_{2}}(t) \Phi_{q_1}(\omega_i) \Phi_{r}(\omega_i) \Phi_{q_2}(\omega_j) \Phi_{s}(\omega_j)\\
\end{aligned}
\end{equation}
Taking the sum of the the bispectrum over the frequency domain yields
\begin{equation}
\begin{aligned}
    & \sum_{k = 1}^{N_\omega} \sum_{i+j=k}^{i \geq 0, j \geq 0} B(t, \omega_i, \omega_j) = \sum_{k = 1}^{N_\omega} \sum_{i+j=k}^{i \geq 0, j \geq 0} \sum_{r=1, s=1}^{N_q, N_q} \sum_{q_1 = 1}^{N_q} \sum_{q_2 = 1}^{N_q} b_{rs}(t) a_{q_{1}}(t) a_{q_{2}}(t) \Phi_{q_1}(\omega_i) \Phi_{r}(\omega_i) \Phi_{q_2}(\omega_j) \Phi_{s}(\omega_j) \\
    & \sum_{k = 1}^{N_\omega} \sum_{i+j=k}^{i \geq 0, j \geq 0} B(t, \omega_i, \omega_j) = \sum_{r=1, s=1}^{N_q, N_q} \sum_{q_1 = 1}^{N_q} \sum_{q_2 = 1}^{N_q} b_{rs}(t) a_{q_{1}}(t) a_{q_{2}}(t) \sum_{k = 1}^{N_\omega} \sum_{i+j=k}^{i \geq 0, j \geq 0} \Phi_{q_1}(\omega_i) \Phi_{r}(\omega_i) \Phi_{q_2}(\omega_j) \Phi_{s}(\omega_j) \\
\end{aligned}
\end{equation}
The orthogonality conditions in this scenario are
\begin{equation}
    \sum_{k = 1}^{N_\omega} \sum_{i + j = k}^{i \geq 0, j \geq 0} \Phi_{r}(\omega_i)\Phi_{q_1}(\omega_i)\Phi_{s}(\omega_j)\Phi_{q_2}(\omega_j) \simeq \begin{cases}
    1, \text{if } q_1 = r, q_2 = s\\
    0, \text{otherwise}
    \end{cases}
\end{equation}
Substituting the orthogonality conditions into the above equation yields
\begin{equation}
    \sum_{k = 1}^{N_\omega} \sum_{i+j=k}^{i \geq 0, j \geq 0} B(t, \omega_i, \omega_j) = \sum_{r=1, s=1}^{N_q, N_q} b_{rs}(t) a_{r}(t) a_{s}(t)
\end{equation}
The definition of the third moment of the process is
\begin{equation}
\begin{aligned}
    E[X^{3}(t)] & = 6 \sum_{k = 1}^{N_\omega} \sum_{i+j=k}^{i \geq 0, j \geq 0} \Re(B(t, \omega_i, \omega_j)) \Delta \omega^{2} \approx 6 \sum_{r=1, s=1}^{N_q, N_q} \Re(b_{rs}(t)) a_{r}(t) a_{s}(t) \Delta\omega^{2} \\
    & = 6 \sum_{r=1, s=1}^{N_q, N_q} |b_{rs}(t)| a_{r}(t) a_{s}(t) \cos(\gamma_{rs}(t)) \Delta\omega^{2} \\
    \label{eqn:third_moment_pod_theoretical}
\end{aligned}
\end{equation}

\section{Ensemble Properties of POD formula}
\label{ap:E}

Here, we show that the expansion Eq.\ \eqref{eqn:pod_expansion} satisfies the necessary ensemble properties in 2nd-order (variance) and 3rd-order respectively as presented in \ref{ap:D}.

\subsection{2-point Correlation Function}

The ensemble 2nd-order properties of the pure and interactive components of the  (expansion are investigated separately.

\subsubsection{Pure Component}

The variance of the pure component of the process can be expanded as
\begin{equation}
    Var_{p}[X(t)] = 4\mathbb{E}[\Delta \omega \sum_{k_1}\sum_{r_1}\sum_{q_2}\sum_{r_2} a_{r_1}(t)a_{r_2}(t)\Phi_{r_1}(\omega_{k_1}) \Phi_{r_2}(\omega_{k_2}) \Delta \omega\cos(\omega_{k_1}t - \phi_{r_1k_1}) \cos(\omega_{k_2}t - \phi_{r_2k_2})]
\end{equation}
and we know that
\begin{equation}
    \mathbb{E}[\cos(\omega_{k_1}t - \phi_{r_1k_1}) \cos(\omega_{k_2}t - \phi_{r_2k_2})] = \begin{cases}
    \frac{1}{2}, \text{if } r_1 = r_2, k_1 = k_2\\
    0, \text{otherwise}
    \end{cases}
\end{equation}
Substituting these values into the ensemble equation yields
\begin{equation}
\begin{aligned}
    Var_{p}[X(t)] & = 2 \Delta \omega \sum_{k}\sum_{r} a_{r}(t)a_{r}(t)\Phi_{r}(\omega_{k}) \Phi_{r}(\omega_{k}) \\
    & = 2 \Delta \omega \sum_{k}\sum_{r} a_{r}^{2}(t) \Phi_{r}^{2}(\omega_{k}) = 2 \Delta \omega \sum_{r} a_{r}^{2}(t) \\
\end{aligned}
\end{equation}
which is the equal to the theoretical value Eq.\ \eqref{eqn:pure_variance_pod_theoretical}.

\subsubsection{Interactive Component}

The variance of the interactive component of the process can be expanded as
\begin{equation}
\begin{aligned}
    & Var_{I}[X(t)] = 4 \mathbb{E}[\sum_{k_1}\sum_{k_2}\sum_{r_1}\sum_{s_1}\sum_{r_2}\sum_{s_2} \sum_{i_1 + j_1 = k_1}^{i_1 \geq j_1 \geq 0} \sum_{i_2 + j_2 = k_2}^{i_2 \geq j_2 \geq 0} |b_{r_1s_1}(t)| |b_{r_2s_2}(t)| \\
    & \Phi_{r_1}(\omega_{i_1}) \Phi_{r_2}(\omega_{i_2}) \Phi_{s_1}(\omega_{j_1}) \Phi_{s_2}(\omega_{j_2}) \Delta \omega^2 \cos(\omega_{k_1}t - \phi_{r_1i_1} - \phi_{s_1j_1} + \gamma_{r_{1}s_{1}}(t)) \cos(\omega_{k_1}t - \phi_{r_2i_2} - \phi_{s_2j_2} + \gamma_{r_{2}s_{2}}(t))]\\
\end{aligned}
\end{equation}
and we know that
\begin{equation}
\begin{aligned}
    & \mathbb{E}[\cos(\omega_{k_1}t - \phi_{r_1i_1} - \phi_{s_1j_1} + + \gamma_{r_{1}s_{1}}(t)) \cos(\omega_{k_1}t - \phi_{r_2i_2} - \phi_{s_2j_2} + + \gamma_{r_{2}s_{2}}(t))] = \\
    & \begin{cases}
    \frac{1}{2}, \text{if } r_1 = r_2, s_1 = s_2, i_1 = i_2, j_1 = j_2, k_1 = k_2\\
    0, \text{otherwise}
    \end{cases} \\
\end{aligned}
\end{equation}
Substituting the values in the ensemble equation yields
\begin{equation}
\begin{aligned}
    Var_{I}[X(t)] & = 2 \sum_{k} \sum_{r}\sum_{s} \sum_{i + j = k}^{i \geq j \geq 0} |b_{rs}(t)|^{2} \Phi_r^2(\omega_i) \Phi_s^2(\omega_j) \Delta \omega^2 \\
    & = \sum_{r}\sum_{s} b_{rs}^2(t) \Delta \omega^2 \sum_{k} \sum_{i + j = k}^{i \geq 0, j \geq 0} \Phi_r^2(\omega_i) \Phi_s^2(\omega_j) \\
    & = \sum_{r}\sum_{s} |b_{rs}(t)|^{2} \Delta \omega^2 \\
\end{aligned}
\end{equation}
which is the equal to the theoretical value Eq.\ \eqref{eqn:interactive_variance_pod_theoretical}

\subsection{Third Moment}
The third moment of the process can be expanded as
\begin{equation}
\begin{aligned}
    \mathbb{E}[X^{3}(t)] & = 48 \mathbb{E}[\sum_{k_1} \sum_{k_2} \sum_{k_3} \sum_{p} \sum_{q} \sum_{r} \sum_{s} \sum_{i_3 + j_3 = k_3}^{i_3 \geq j_3 \geq 0} a_{p}(t) \Phi_{p}(\omega_{k_{1}}) a_{q}(t) \Phi_{q}(\omega_{k_{2}}) |b_{rs}(t)| \\
    & \theta_{rs}(\omega_{i_3}, \omega_{j_3}) \Delta \omega^2 \cos(\omega_{k_1}t - \phi_{pk_1}) \cos(\omega_{k_2}t - \phi_{qk_2}) \cos(\omega_{k_3}t - \phi_{ri_3} - \phi_{sj_3} + \gamma_{rs}(t))] \\
\end{aligned}
\end{equation}
and we know that
\begin{equation}
\begin{aligned}
    & \mathbb{E}[\cos(\omega_{k_1}t - \phi_{pk_1}) \cos(\omega_{k_2}t - \phi_{qk_2}) \cos(\omega_{k_3}t - \phi_{ri_3} - \phi_{sj_3} + \gamma_{rs}(t))] = \\
    & \begin{cases}
    \frac{1}{4} \cos(\gamma_{rs}(t)), \text{if } p = r, q = s, k_1 = i_3, k_2 = j_3\\
    0, \text{otherwise}
    \end{cases} \\
\end{aligned}
\end{equation}
Substituting this expectation into the ensemble equation yields
\begin{equation}
\begin{aligned}
    \mathbb{E}[X^{3}(t)] & = 12 \sum_{k} \sum_{r} \sum_{s} \sum_{i + j = k}^{i \geq j \geq 0} a_{r}(t) \Phi_{r}(\omega_i) a_{s}(t) \Phi_{s}(\omega_j) |b_{rs}(t)| \cos(\gamma_{rs}(t)) \Phi_{r}(\omega_{i}) \Phi_{s}(\omega_{j}) \Delta \omega^2 \\
    & = 6 \sum_{r} \sum_{s} a_{r}(t) a_{s}(t) |b_{rs}(t)| \cos(\gamma_{rs}(t)) \Delta \omega^2 \\
\end{aligned}
\end{equation}
which is the same as Eq.\ \eqref{eqn:third_moment_pod_theoretical}

\end{document}